\documentclass[ijoc,nonblindrev]{informs3} 

\OneAndAHalfSpacedXII 

\usepackage{endnotes}
\usepackage{comment}
\usepackage{subfigure}

\let\footnote=\endnote

\newcommand{\eqdef}{\stackrel{\rm def}{=}}

\usepackage{natbib}
 \bibpunct[, ]{(}{)}{,}{a}{}{,}%
 \def\bibfont{\small}%
 \def\newblock{\ }%
\EquationsNumberedThrough    

\newtheorem{theorem}{Theorem}[section]
\newtheorem{lemma}{Lemma}[section]

\newtheorem{assumption}{Assumption}[section]
\newtheorem{proposition}{Proposition}[section]

\begin{document}



\RUNTITLE{Approximate Dynamic Programming Methods for Energy Storage}

\TITLE{         \Large{Least Squares Policy Iteration with Instrumental Variables }\\
\vspace{0.00 in} \Large{vs. Direct Policy Search: Comparison Against Optimal }\\
\vspace{0.00 in} \Large{Benchmarks Using Energy Storage}}




\ARTICLEAUTHORS{%
\AUTHOR{Warren R. Scott}
\AFF{Department of Operations Research and Financial Engineering, Princeton University, Princeton, NJ 08544, \EMAIL{wscott@princeton.edu}}
\AUTHOR{Warren B. Powell, Somayeh Moazeni}
\AFF{Department of Operations Research and Financial Engineering, Princeton University, Princeton, NJ 08544, \EMAIL{powell@princeton.edu, somayeh@princeton.edu}}
} 

\ABSTRACT{This paper studies approximate policy iteration (API) methods which use least-squares Bellman error minimization for policy evaluation. We address several of its enhancements, namely, Bellman error minimization using instrumental variables, least-squares projected Bellman error minimization, and projected Bellman error minimization using instrumental variables. We prove that for a general discrete-time stochastic control problem, Bellman error minimization using instrumental variables is equivalent to both variants of projected Bellman error minimization. An alternative to these API methods is direct policy search based on knowledge gradient. The practical performance of these three approximate dynamic programming methods are then investigated in the context of an application in energy storage, integrated with an intermittent wind energy supply to fully serve a stochastic time-varying electricity demand. We create a library of test problems using real-world data and apply value iteration to find their optimal policies. These benchmarks are then used to compare the developed policies. Our analysis indicates that API with instrumental variables Bellman error minimization prominently outperforms API with least-squares Bellman error minimization. However, these approaches underperform our direct policy search implementation.
}


\HISTORY{}

\maketitle

%

\section{Introduction}

We have long recognized that the powerful theory of discrete Markov decision processes \citep{Pu94} is limited by the well known curse of dimensionality, which can be traced to the need to compute the value of being in each discrete state or, more problematically, the probability of transitioning from one discrete state to another given an action.  Known as a {\it flat representation} in computer science, even small problems quickly blow up using this model.  The field of approximate dynamic programming (\cite{BeTs96}, \cite{Po11}, \cite{BertsekasADP2012}) and reinforcement learning (\cite{SuBa98}, \cite{Szepesvari2010b}) have offered the hope of partially overcoming this problem by replacing the value function using a statistical approximation (in particular, linear regression), allowing us to draw on approximate versions of powerful algorithmic strategies such as value iteration and policy iteration \citep{Pu94}.

Somewhat surprisingly, while we have found in our own work that approximate value iteration works well for very specific problem classes (\cite{ToPo06}, \cite{PoGeLa11}, \cite{SiDaGe09}, \cite{NaPo2013}), it does not work as a general algorithmic strategy (our success with AVI has always been in the context of problems where we could exploit convexity).  The computer science community has focused considerable efforts on a strategy called $Q$-learning, which involves learning the value of a state-action pair rather than just the value of being in a state (\cite{SuBa98}).  $Q$-learning enjoys rigorous convergence theory for lookup table representations \citep{Ts94}, but this does not scale even to small problems.  Approximating $Q$-factors with linear models is more art than science. Convergence results (\cite{SuSzMa09}, \cite{Maei2009}) do not contain any performance guarantees (they only establish convergence of a particular approximating architecture), and we are unaware of any empirical comparisons against optimal benchmarks.


Perhaps for this reason, approximate policy iteration has attracted considerable recent attention (see \cite{BertsekasADP2012} and the review in \cite{PoMa11}). API avoids the need to approximate the value of a state-action pair, and enjoys stronger convergence theory, although this always involves assumptions that are unlikely to be perfectly satisfied in practice.  However, we are still unaware of any comparisons against optimal benchmarks.  In this paper, we use the setting of a class of energy storage problems that requires balancing power from the grid and power from a stochastic, renewable source, to serve a load (that is sometimes time varying), with access to an energy storage device.  While the datasets were carefully generated, the focus of the paper is purely algorithmic.  Energy storage is a problem that is well suited to the use of value function approximations.  Further, while the problem is realistic, we could, with minor simplifications, create problems that could be solved optimally, providing us with rigorous benchmarks that allow us to provide an accurate assessment of the quality of the solution produced by different algorithmic strategies.  The software containing these benchmark problems is available at http://www.castlelab.princeton.edu/datasets.htm.



In this paper, we focus on the use of linear architectures for approximating the value function, a strategy that has received the most attention in the literature.  We model our problem in steady state, which allows us to use a powerful algorithmic strategy called Least Squares Approximate Policy Iteration (LSAPI), introduced by \cite{lagoudakis2003}, building on the least squares temporal difference method, introduced by \cite{BrBa96} to estimate the value of a fixed policy.  \cite{lagoudakis2003} uses a statistical technique that introduces the idea of using instrumental variables (IV) to overcome errors in the observation of the independent variables (explained below).


Recent research has revealed the importance of using the idea of modifying Bellman's optimality equation by first projecting the value function onto the space of linear architectures.  This is referred to in the literature as {\it Mean-Squared Projected Bellman Error} (MSPBE) minimization (see \cite{LaPa03}, \cite{SuMaei2009}, \cite{MaSz10}, and the explanation of this material in \cite{Po11}[Chapter 9]).  In this paper, we show for the first time that MSPBE is mathematically equivalent to the use of instrumental variables, which is also simpler.

Although our focus is on using exact benchmarks to derive insights into different algorithmic strategies, our choice of energy storage is motivated by the importance of this problem class.  For example, \cite{So81} addresses the potential benefits of combining a hydro storage facility with wind energy. Three control policies and their corresponding (expected) gained revenues are analyzed in \citep{Storage:IntermittentEnergy:2004}, assuming probabilistic models for load and wind, and a load-price curve for the electricity price. \cite{CaLu05} and \cite{LaMaSe10} investigate the operation and value of a gas storage device. \cite{Gr07} and \cite{Sw07} discuss incorporating a compressed air energy storage with wind power generation in energy systems. Maximization of the expected market profit of a wind farm and a hydro pumped storage owner over a finite horizon to comply with his commitments in the market is addressed in \citep{GMSG2008}. The problem is formulated as a two-stage stochastic programming problem with random market prices and random wind generation. \cite{Si09} studies the potential value of a storage device in the PJM network used for arbitrage. For a thorough review of this growing literature, we refer the reader to \citep{Harsha:2012} and the references therein.  The paper that addresses a problem closest to ours is \cite{sioshansi2014}, but this paper is able to use classical backward dynamic programming using lookup tables, and does not explore approximation strategies that might scale to more complex problems.  

For our energy storage application, we created a library of test problems using realistic data, constrained by the goal of creating optimal benchmarks. Optimal policies for these problems are computed using exact value iteration \citep{Pu94}. These benchmarks are then used to compare variations policies estimated using approximate policy iteration and direct policy search. Our analysis indicates that API with Bellman error minimization using instrumental variables significantly outperforms API with least-squares Bellman error minimization.



This paper makes the following contributions.  1) We introduce and prove the consistency of a hybrid policy based on least squares approximate policy iteration that combines instrumental variables and projected Bellman error minimization. 2) We establish, for the first time, the equivalence of least-squares approximate policy iteration using instrumental variables with mean-squared projected Bellman error minimization, as well as the new hybrid policy combining both instrumental variables and projected Bellman error minimization.  3) We develop a set of benchmark problems using the important context of energy storage where we are able to derive optimal policies.  While these problems are relatively simple, CPU times to estimate these policies typically ranged around two {\it weeks}.  4) We calculated approximate policies using basic LSAPI, LSAPI using instrumental variables, and direct policy search, and compared the results against optimal benchmark policies.  These experiments show that LSAPI with instrumental variables performs significantly better than basic LSAPI, and yet even this advanced algorithmic strategy falls far short of optimal.  Direct policy search outperforms both LSAPI policies by a wide margin.  Since the structure of the policy is the same, the issue is not the accuracy of the approximating architecture, but rather the estimation procedure, calling into question the validity of Bellman error minimization.

This paper is organized as follows. Section \ref{sec:APIA} provides an overview of approximate policy iteration methods. Several policy evaluation algorithms based on Bellman error minimization are discussed in Section \ref{Sec:LSTD:PE}. The proof for their equivalence is presented in Section \ref{sec:PolicyEvaluationAlgsEquivalent}. Direct policy search and knowledge gradient are discussed in Section \ref{sec:DirectPolicySearch}. The energy storage management problem, its underlying stochastic processes, and its stochastic dynamic programming formulation are explained in Section \ref{sec:Models}. Performance of the computed approximate dynamic programming policies is investigated and compared with benchmark problems in Section \ref{sec:NumericalWork}. Concluding remarks are given in Section \ref{Sec:conclusions}.

\section{Approximate Policy Iteration Algorithms}\label{sec:APIA}
Stochastic dynamic programming for maximizing expected revenues relies on Bellman's optimality equation given by
\begin{eqnarray}
\label{eq:BellmansEquationPre}&&V(S_t) = \displaystyle\max_{x}~ \mathbb{E} \left[ C\left(S_t,x\right) + \gamma~ V\left(S_{t+1}\right)|S_{t}\right],
\end{eqnarray}
where $S_{t}$ refers to the state variable at time step $t$, $V(\cdot)$ is the optimal value function (around the pre-decision state $S_t$), $C(S_{t},x)$ is the contribution function, and $0\leq \gamma< 1$ is the discount factor. The expectation in equation \eqref{eq:BellmansEquationPre} is over the exogenous random changes, denoted by $W_{t+1}$, in the state of the system. The state variable at the next step is then obtained by the transition function $S_{t+1} = S^{M}\left(S_t,x,W_{t+1}\right)$. Throughout, we use the convention that any variable indexed by $t$ is known at time $t$.

Often computing the expected value of $V(S_{t})$ is challenging and needs to be approximated. To avoid this, a modified version of Bellman's equation based on the \emph{post-decision state variables} can be adopted, see e.g., \citep{Be11} or Section $4.6$ of \citep{Po11} for a thorough discussion of post-decision states. The post-decision state variable, denoted by $S^{x}_{t}$, refers to the state immediately after being in the pre-decision state $S_t$ and taking the action/decision $x$, but before any exogenous information (randomness) from the state transition has been revealed. The \emph{post-decision value function}, denoted by $V^{x}(S^{x}_{t})$, explains the value of being in the post-decision state $S^x_{t}$. The relationship between $V^{x}(S^{x}_{t})$ and $V(S_{t})$ is defined as:
\begin{eqnarray}
\nonumber&&V^x\left(S^x_{t}\right) \eqdef \mathbb{E}\left[V\left(S_{t+1}\right)~|~S^x_{t}\right].
\end{eqnarray}
Thus, Bellman's equation \eqref{eq:BellmansEquationPre} around the post-decision state variables can be written as
\begin{eqnarray}
\label{eq:BellmansEquationPost}&&V^{x}_{t-1}(S^{x}_{t-1}) = \mathbb{E}\left[\displaystyle\max_x ~ \left( C\left(S_t, x\right) + \gamma~ V^x_{t}\left(S^x_{t}\right) \right) ~|~ S^x_{t-1}\right].
\end{eqnarray}
The expectation being outside of the maximum operator allows us to solve the inner maximization problem using deterministic optimization techniques.   Furthermore, in some problems, the post-decision state has a lower dimension than the pre-decision state.

For most applications, however, when the state variable is multidimensional and continuous, Bellman's equations \eqref{eq:BellmansEquationPre} or \eqref{eq:BellmansEquationPost} cannot be solved exactly; as a result a large field of research on approximation techniques has evolved, see e.g., \citep{Neuro:Dynamic,SuBa98,Szepesvari2010b,Po11}. In this paper, we focus on the widely used strategy of linear architectures of the form
\begin{eqnarray}
\label{EQ:ApproxValFunc}&&\overline{V}^{x}_{t}(S^x_{t})\eqdef \sum_{f\in\mathcal{F}} \theta_{f}~\phi_{f}\left(S^x_{t}\right)= \theta^{\top} \phi\left(S^x_{t}\right),
\end{eqnarray}
where $ \phi\left(S^x_{t}\right)$ is the column vector with elements $\left\{\phi_{f}(S^{x}_{t})\right\}_{f\in\mathcal{F}}$, and $\theta$ is the column vector of weights for the basis functions. Substituting approximate post-decision value functions \eqref{EQ:ApproxValFunc} into equation \eqref{eq:BellmansEquationPost}, for a fixed policy $\pi$ and a given vector $\theta$, we get
\begin{eqnarray}
\label{eq:BellmansEquationPostParametric}&&\theta^{\top}\phi\left(S^x_{t-1}\right) = \mathbb{E}\left[ C\left(S_t, X^\pi\left(S_t | \theta\right)\right) + \gamma~\theta^{\top}\phi\left(S^x_t\right) ~|~ S^x_{t-1}\right].
\end{eqnarray}
In general, a fixed point satisfying this equation does not exist using linear architectures \citep{FaVr00}, but this equation forms the foundation of LSAPI.  For a given weight vector $\theta$, our policy $X^\pi(S_t | \theta)$ is given by
\begin{eqnarray}
\label{X:Def}&& X^\pi(S_t|\theta) = \text{argmax}_x~ \left[C\left(S_t,x\right) +  \gamma~\theta^{\top}\phi\left(S^x_t\right)\right].
\end{eqnarray}


\begin{figure}
\begin{center}
\begin{tabular}{|l|}
\hline
Approximate Policy Iteration Algorithm\\
\hline
(01) Initialize $\theta$.\\
(02) \textbf{for} $j = 1$ \textbf{to} $M$ \textbf{(Policy Improvement Loop)}\\
(03) \hspace{.5 cm}\textbf{for} $i =1$ \textbf{to} $N$ \textbf{(Policy Evaluation Loop)}\\
(04) \hspace{1.0 cm}Simulate a random post-decision state, $S^{x}_{t-1,i}$.\\
(05) \hspace{1.0 cm}Record $\phi(S^x_{t-1,i})$.\\
(06) \hspace{1.0 cm}Simulate the state transition to get $S_{t,i}$.\\
(07) \hspace{1.0 cm}Determine the decision, $x= X^\pi(S_{t,i}|\theta)$ using equation \eqref{X:Def}.\\
(08) \hspace{1.0 cm}Record $C(S_{t,i}, x).$\\
(09) \hspace{1.0 cm}Record $\phi(S^x_{t,i})$, the observation of $\mathbb{E}[\phi(S^x_{t,i})| S^x_{t-1,i}]$.\\
(10) \hspace{.5cm} \textbf{End}\\
(11) \hspace{.5cm}Update $\theta$ using equation \eqref{eq:BellmanErrorLS}, or \eqref{eq:LSAPImatrix}, or \eqref{eq:LSAPIprojected}, or \eqref{eq:IVProjectedBellmanErrorMinimization1}. \textbf{(Policy Evaluation)}\\
(12) \textbf{End}\\
\hline
\end{tabular}
\end{center}
\caption{Summary of (steady-state) approximate policy iteration.  Policy Evaluation loop approximately evaluates the cost function of
a given policy. Policy Improvement loop generates a new (improved) policy. The constants $M$ and $N$ indicate the maximum number of policy improvement and policy evaluation iterations.}
\label{fig:LSAPI}
\end{figure}

The basic LSAPI algorithm of \cite{BrBa96} is given in figure \ref{fig:LSAPI}.  We use the adaption for the post-decision state given in \cite{MaPo10}. Let $\{ S^x_{t-1,i} \}_{i=1}^N$ be a set of (randomly generated) post-decision states. For each sample $i=1,\cdots,N$, the next post-decision state $S^x_{t,i}$ is simulated and, using equation \eqref{X:Def}, the the optimal decision $X^\pi(S_{t,i}|\theta)$ is computed. For every $i=1,\cdots,N$, equation \eqref{eq:BellmansEquationPostParametric} can be rewritten as:
\begin{eqnarray}
\label{eq:BellmansEquationPostParametricRearranged}&C\left(S_{t,i}, X^\pi(S_{t,i} | \theta)\right)=& \left(\phi(S^x_{t-1,i}) -\gamma \mathbb{E}[\phi(S^x_{t,i})| S^x_{t-1,i}] \right)^{\top} \theta \\
\nonumber&&+ C(S_{t,i}, X^\pi(S_{t,i} | \theta)) - \mathbb{E}[ C(S_{t,i}, X^\pi(S_{t,i} | \theta))| S^x_{t-1,i}].
\end{eqnarray}
Denote $C_{t,i}\eqdef C\left(S_{t,i}, X^\pi(S_{t,i} | \theta)\right)$ and $X_{t,i}\eqdef \left(\phi(S^x_{t-1,i}) -\gamma \mathbb{E}[\phi(S^x_{t,i})| S^x_{t-1,i}] \right)^{\top}$. Using simulation, we are able to get observations of $C_{t,i}$ and $X_{t,i}$ in equation \eqref{eq:BellmansEquationPostParametricRearranged}.

Let $\Phi_{t}$ and $\Phi_{t-1}$ be the matrices of fixed basis functions whose rows correspond to states and columns correspond to basis functions. Furthermore, refer to the vector of contribution functions and its expected value by $C_{t}$ and $\overline{C}_{t}$, respectively. More precisely,
\begin{eqnarray}
\nonumber \Phi_{t-1} \eqdef \begin{bmatrix} \phi(S^x_{t-1,1})^{\top} \\ \vdots \\ \phi(S^x_{t-1,N})^{\top} \end{bmatrix},~~
{\Phi}_{t} \eqdef \begin{bmatrix} \phi(S^x_{t,1})^{\top} \\ \vdots \\ \phi(S^x_{t,N})^{\top} \end{bmatrix},~~
C_{t} \eqdef \begin{bmatrix}C_{t,i}\\ \vdots \\C_{t,N}\end{bmatrix},~~
\overline{C}_{t} \eqdef \begin{bmatrix} \overline{C}_{t,1} \\ \vdots \\ \overline{C}_{t,N}\end{bmatrix},
\end{eqnarray}
where $\overline{C}_{t,i}\eqdef \mathbb{E}[ C(S_{t,i}, X^\pi(S_{t,i} | \theta))| S^x_{t-1,i}]$. Thus, equation \eqref{eq:BellmansEquationPostParametricRearranged} can be presented in a matrix form as
\begin{equation}
\underbrace{C_{t}}_{n \times 1} = \underbrace{(\Phi_{t-1} - \gamma {\Phi}_{t})}_{n \times k} \underbrace{\theta}_{k \times 1} + \underbrace{C_{t}-\overline{C}_{t}}_{n \times 1}, \label{eq:BellmansEquationPostParametricMat}
\end{equation}
We refer to $(C_{t}-\overline{C}_{t})$ as the \emph{Bellman errors}, although this is not the standard definition of the Bellman error.

An improved value for $\theta$ can be obtained by minimizing the Euclidean norm of the Bellman errors,
\begin{eqnarray}
\label{Eq:MainMinError}&&\min_{\theta}~\left\|C_{t}- \overline{C}_{t}\right\|_{2}^{2},
\end{eqnarray}
which is the goal in least-squares Bellman error minimization. In the next section, we address a few of its variants including Bellman error minimization using instrumental variables and projected Bellman error minimization, after which we prove that these two methods are equivalent. In our discussion, we make the following assumption about the linear independence of the basis functions and underlying matrices.
\begin{assumption}
The matrices $\Phi_{t-1},~(\Phi_{t-1} - \gamma {\Phi}_{t})$, and $\Phi_{t-1}^{\top} (\Phi_{t-1} - \gamma {\Phi}_{t})$ have full column rank, and $k \leq n$.\label{ass:PhiFullRank}
\end{assumption}
This assumption ensures that enough different states are visited that allows us to identify the linear model.

\section{Least Squares Temporal Differences (LSTD) Policy Evaluation}\label{Sec:LSTD:PE}
The LSTD method, introduced by \cite{BrBa96}, is one of the (batch) variants of the TD algorithms. For a given policy and assuming a linear model with fixed basis functions for the value function, this on-policy algorithm performs a least squares regression so that the sum of the temporal differences over the simulation to approximate the expectation equals zero. LSTD tends to be statistically efficient and extract more information from training experiences, compared to some other TD methods. The method of instrumental variables, see e.g., \cite{IV:Methods}, is used in LSTD to compensate for the use of Monte Carlo simulation to approximate the expected cost in $\overline{C}_{t,i}\eqdef \mathbb{E}[ C(S_{t,i}, X^\pi(S_{t,i} | \theta))| S^x_{t-1,i}]$. The estimation of the value function converges with probability one to the true value, when the true value function lies in the span of the basis functions, the number of state transitions increases to infinity, and a few other technical assumptions are fulfilled.

\cite{LaPa03} expands upon the LSTD algorithm in \citep{BrBa96} by using a linear architecture to approximate the value function over the higher dimension state-action pairs.  Furthermore, they give a geometric interpretation of several different methods of approximately solving Bellman's equation. Also assuming finite states and actions, \cite{LaPa03} introduces least-squares policy iteration which approximates the value of state-action pairs (Q-factors) with a linear model. When doing policy evaluation, they choose to use least-squares to minimize the fixed-point approximation error instead of the Bellman error. An alternative approach is to project the Bellman errors down onto the space spanned by the basis functions of the value function and then minimize the Bellman errors.  \cite{LaPa03} explains that in general the approximate value function is a fixed point of the projected Bellman operator, not the Bellman operator (see \cite{DeVa00} for a nice discussion).

\cite{BrBa96} addresses conditions for on-policy policy evaluation subproblem based on Bellman error minimization to converge to a fixed value function when using a linear model for the value function. Given a fixed policy, \cite{TsVr97} provide conditions under which on-policy temporal-difference learning will converge to the true value function projected onto the space of value function approximations (with respect to a weighted norm). Furthermore, \cite{TsVr97} proves that, when evaluating a fixed policy, off-policy TD($\lambda$) algorithms with a linear function approximation of the value function may not converge. In this case, off-policy means that the distribution of the states visited during a single infinite trajectory is not equal to the distribution of the states visited if we followed the fixed policy of interest. However, \cite{LaPa03} explain that on-policy LSTD biases the value function and may fail to acceptably fit the value function at states that are rarely visited. Another major disadvantage of on-policy learning is that Assumption \ref{ass:PhiFullRank} may not hold if the policy does not explore enough states.  An important point to keep in mind is that the value of the greedy policy determined by the final value function approximation may be significantly different from the true value function.

In the following, we summarize the policy evaluation algorithms based on Bellman error minimization and instrumental variables. We focus on the off-policy case in both subsections \eqref{sec:BellmanErrorMinimization} and \eqref{sec:ProjectedBellmanErrorMinimization}. Subsection \ref{sec:BellmanErrorMinimization} addresses least-squares Bellman error minimization and instrumental variables Bellman error minimization. Subsection \ref{sec:ProjectedBellmanErrorMinimization} explains least-squares projected Bellman error minimization and instrumental variables projected Bellman error minimization.

\subsection{Policy Evaluation using Bellman Error Minimization}\label{sec:BellmanErrorMinimization}
Applying the typical method of least-squares, a solution of equation \eqref{Eq:MainMinError} equals,
\begin{eqnarray}
\label{eq:BellmanErrorLS}&&\hat{\theta} = \left[(\Phi_{t-1} - \gamma {\Phi}_{t})^{\top} (\Phi_{t-1} - \gamma {\Phi}_{t})\right]^{-1} (\Phi_{t-1} - \gamma {\Phi}_{t})^{\top} C_{t},
\end{eqnarray}
which we refer to as the \emph{least-squares Bellman error minimization} estimator. The matrix of regressors, $(\Phi_{t-1} - \gamma {\Phi}_{t})$, is not deterministic (${\Phi}_{t}$ is not deterministic because we cannot calculate $\mathbb{E}[\phi(S^x_t)| S^x_{t-1}]$); we can only simulate $\phi(S^x_t)$ given $S^x_{t-1}$ and, as a result, the least-squares estimator for $\theta$ will typically be inconsistent.

A class of simple and computationally efficient techniques to obtain consistent estimates, without modeling the noise, is referred to as \emph{Instrumental Variable Methods} (IVM). An instrumental variable is a variable that is correlated with the regressors, but uncorrelated with the errors in the regressors and the observations, see e.g., \citep{Du54,KeSt61,So83,BoTu90}. We provide a brief overview of IV methods in Appendix \ref{app:IV}. Instrumental variables have been previously used in the context of approximate policy iteration algorithms, see e.g., \citep{BrBa96,MaPo10}, but otherwise have not received much attention, even in the reinforcement learning literature.

This results in what we call \emph{instrumental variables Bellman error minimization},
\begin{eqnarray}
\label{eq:LSAPImatrix}\hat{\theta} = \left[(\Phi_{t-1})^{\top} (\Phi_{t-1} - \gamma {\Phi}_{t})\right]^{-1}(\Phi_{t-1})^{\top} C_{t}.
\end{eqnarray}
\cite{BrBa96} provide conditions such that equation \eqref{eq:LSAPImatrix} is a consistent estimator ($\lim_{n \to \infty} \hat{\theta} = \theta$ with probability one) for the on-policy case.  The proof references the consistency properties of the method of instrumental variables by showing that the columns of $\Phi^n$ are appropriate instrumental variables (see Appendix \ref{app:IV}).  One interesting comment is that the matrix $[\Phi_{t-1}^{\top} (\Phi_{t-1} - \gamma {\Phi}_{t})]$ could have negative eigenvalues, unlike $\left[(\Phi_{t-1} - \gamma {\Phi}_{t})^{\top} (\Phi_{t-1} - \gamma {\Phi}_{t})\right]$.

\subsection{Policy Evaluation using Projected Bellman Error Minimization} \label{sec:ProjectedBellmanErrorMinimization}
The idea of projected Bellman error minimization, also called \emph{least-squares fixed-point approximation} in \cite{LaPa03}, is to first project the Bellman errors onto the space spanned by the basis functions defining the value function and then minimize them, see \citep{LaPa03,SuMa09}.  Projecting the left and right hand sides of equation \eqref{eq:BellmansEquationPostParametricMat} onto the space spanned by $\Phi_{t-1}$ with respect to the Euclidean norm, we get
\begin{equation}
\Pi_{t-1} C_{t} = \Pi_{t-1} (\Phi_{t-1} - \gamma \Phi_{t}) \theta + \Pi_{t-1} (C_{t}-\overline{C}_{t}), \label{eq:BellmansEquationPostParametricVecProj}
\end{equation}
where $\Pi_{t-1} = \Phi_{t-1} ((\Phi_{t-1})^{\top} \Phi_{t-1})^{-1} (\Phi_{t-1})^{\top}$ is the projection operator onto the space spanned by the basis functions (see \citep{TsVr97} for the original derivation of this mapping or Section 8.2.3 of \cite{Po11}). Assumption \ref{ass:PhiFullRank} implies that the matrix $\Phi_{t-1}$ has full column rank, and hence $\Pi_{t-1}$ is well defined. We refer to $\Pi_{t-1} (C_{t}-\overline{C}_{t})$, as the \emph{projected Bellman error}. Notice that the $n\times n$ matrix $\Pi_{t-1}$ is the projection with respect to the Euclidean norm that solves the problem: $\min_{\Pi} \| \Phi_{t-1} \theta - b \|_2 = \| \Pi_{t-1} b - b \|_2$. In other words, to get from an arbitrary vector $b$ to the closest vector (in the Euclidean norm sense) which is in the span of the columns of $\Phi_{t-1}$, it is sufficient to just apply the projection $\Pi_{t-1}$ to the vector $b$.

Typically, equation \eqref{eq:BellmansEquationPostParametricVecProj} is an over-determined set of equations.  Taking a least squares approach, we find $\theta$ by minimizing the norm of the projected Bellman error
\begin{eqnarray}
\min_{\theta}~ \left\| \Pi_{t-1} (C_{t} -\overline{C}_{t}) \right\|_2
&=& \min_{\theta} ~\left\| \Pi_{t-1} C_{t} - \Pi_{t-1} (\Phi_{t-1} - \gamma \Phi_{t}) \theta \right\|_2. \nonumber
\end{eqnarray}
The least-squares estimator of $\theta$ then yields
\begin{eqnarray}
\label{eq:LSAPIprojected}&&\hat{\theta} = \left[\left(\Pi_{t-1} (\Phi_{t-1} - \gamma \Phi_{t})\right)^{\top} \left(\Pi_{t-1} (\Phi_{t-1} - \gamma \Phi_{t})\right)\right]^{-1} \left(\Pi_{t-1} (\Phi_{t-1} - \gamma \Phi_{t})\right)^{\top} \Pi_{t-1} C_{t},
\end{eqnarray}
to which we refer as \emph{least-squares projected Bellman error minimization}.

To establish a consistent estimator for $\theta$, similar to the derivation of equation \eqref{eq:LSAPImatrix}, $\Phi_{t-1}$ can be used as an instrumental variable, e.g., see Appendix \ref{app:IV} or the proof in \citep{BrBa96}. We refer to the resulting estimator as the \emph{instrumental variables projected Bellman error minimization}, 
\begin{eqnarray}
\label{eq:IVProjectedBellmanErrorMinimization1}&&\hat{\theta} = \left[ (\Phi_{t-1})^{\top} \Pi_{t-1} (\Phi_{t-1} - \gamma \Phi_{t}) \right]^{-1} (\Phi_{t-1})^{\top} \Pi_{t-1} C_{t}.
\end{eqnarray}
We note that $\Pi_{t-1}\Phi_{t-1}$ could also have been used for the instrumental variables instead of $\Phi_{t-1}$. However, it is easy to see that the obtained estimator would be equivalent to that in equation \eqref{eq:IVProjectedBellmanErrorMinimization1}.

In the rest of this section, we show that the estimator $\hat{\theta}$ in equation \eqref{eq:IVProjectedBellmanErrorMinimization1} is consistent (converges in probability to the true weights), for the off-policy case. \cite{BrBa96} proves a similar result for the on-policy case. The following discussion remains valid even when the state space is continuous or the discount factor is one.

We first make the following assumptions:
\begin{assumption}\label{ass:OffPolicy}
The rows of $\Phi_{t-1}$ are i.i.d. (off-policy).
\end{assumption}
\begin{assumption}\label{ass:IVsigmaFullRankDP}
The covariance matrix $\Sigma$ between the instrumental variables and regressors has full rank $k$, where $\Sigma_{jl} = \text{Cov}[(\Phi_{t-1})_{j}, ~X_l]$ and $X= \Pi_{t-1} (\Phi_{t-1} - \gamma \mathbb{E}[ \gamma \Phi_{t} \vert \{S^x_{t-1}\}])$.
\end{assumption}
\begin{assumption}\label{ass:IVuncorrelatedWithYnoiseLimitLLN}
$\mathbb{E}[| (\Phi_{t-1})_{ij} Y^{''}_i |] < \infty$, for $j=1,2,\cdots,k$, where $Y'' = \Pi_{t-1} (C_{t}-\overline{C}_{t})$.
\end{assumption}
\begin{assumption}\label{ass:IVuncorrelatedWithXnoiseLimitLLN}
$ \mathbb{E}[| (\Phi_{t-1})_{ij} \left(\Pi_{t-1} (\Phi_{t-1} - \gamma \Phi_{t}) - X\right)_{il}|] < \infty$, for $j,l \in \{1,2,\cdots,k\}$, where the matrix $X$ as in Assumption \eqref{ass:IVsigmaFullRankDP}.
\end{assumption}
Assumptions \ref{ass:IVuncorrelatedWithYnoiseLimitLLN} and \ref{ass:IVuncorrelatedWithXnoiseLimitLLN} are necessary to use the law of large numbers, which guarantees sample means converge to their true means.

\begin{theorem}\label{Theorem:consistency}
Let assumptions \ref{ass:PhiFullRank}, \ref{ass:OffPolicy}, \ref{ass:IVsigmaFullRankDP}, \ref{ass:IVuncorrelatedWithYnoiseLimitLLN}, and \ref{ass:IVuncorrelatedWithXnoiseLimitLLN} hold. The estimator $\hat{\theta}$ presented in equation \eqref{eq:IVProjectedBellmanErrorMinimization1} is a consistent estimator for $\theta$ defined in equation \eqref{eq:BellmansEquationPostParametricVecProj}.

\textbf{Proof:} The proof follows from Theorem \ref{thm:IVconsistent} in Appendix \ref{app:IV}.  The following lemmas show that the assumptions in Appendix \ref{app:IV} for Theorem \ref{thm:IVconsistent} hold. For notation consistent with Appendix \ref{app:IV}, we let $X\eqdef \Pi_{t-1} (\Phi_{t-1} - \gamma \mathbb{E}[ \gamma \Phi_{t} \vert \{S^x_{t-1}\}])$, $X'\eqdef \Pi_{t-1} (\Phi_{t-1} - \gamma \Phi_{t})$, $X'' \eqdef X' - X$, $Y'' \eqdef \Pi_{t-1} (C_{t}-\overline{C}_{t})$, and $Z\eqdef \Phi_{t-1}$. We first show Assumption \ref{ass:IVmeanYnoiseZero} holds.

\begin{lemma}\label{lem:YnoiseMean0}
For every $i$, $\mathbb{E}[Y^{''}_{i}]=0$.
\newline \textbf{Proof:}
\begin{eqnarray}
\nonumber\mathbb{E}[Y^{''}] &=& \mathbb{E}[\Pi_{t-1} (C_{t}-\overline{C}_{t})] \\
\nonumber&=& \mathbb{E}\left[\mathbb{E}[\Pi_{t-1} (C_{t}-\overline{C}_{t}) \vert \{S^x_{t-1}\}] \right]\\
\nonumber&=& \mathbb{E}[ \Pi_{t-1} \underbrace{\mathbb{E}[(C_{t}-\overline{C}_{t}) \vert \{S^x_{t-1}\}]}_{= \vec{0}} ]=\vec{0}.
\end{eqnarray}
\end{lemma}

We next show Assumption \ref{ass:IVmeanXnoiseZero} holds, which states that the mean of the noise in the observation of the explanatory variables is zero.
\begin{lemma}\label{lem:XnoiseMean0}
For every $i,j$, $\mathbb{E}[X^{''}_{ij}]=0$.
\newline \textbf{Proof:}
\begin{eqnarray}
\nonumber\mathbb{E}[X^{''}] &=& \mathbb{E}[X^{'} - X] \\
\nonumber&=& \mathbb{E}\left[\Pi_{t-1} (\Phi_{t-1} - \gamma \Phi_{t}) - \Pi_{t-1} (\Phi_{t-1} - \mathbb{E}[ \gamma \Phi_{t} \vert \{S^x_{t-1}\}])\right]\\
\nonumber &=&  \gamma \mathbb{E}\left[\Pi_{t-1}(\mathbb{E}[\Phi_{t} \vert \{S^x_{t-1}\}] - \Phi_{t})\right]\\
\nonumber&=&  \gamma \mathbb{E}\left[\mathbb{E}[  \Pi_{t-1}(\mathbb{E}[\Phi_{t} \vert \{S^x_{t-1}\}] - \Phi_{t}) \vert \{S^x_{t-1}\} ]\right]\\
\nonumber&=&  \gamma \mathbb{E}\left[ \Pi_{t-1} \mathbb{E}[  \mathbb{E}[\Phi_{t} \vert \{S^x_{t-1}\}] - \Phi_{t} \vert \{S^x_{t-1}\} ]\right]\\
\nonumber&=&  \gamma \mathbb{E}[ \Pi_{t-1} (\underbrace{\mathbb{E}[\Phi_{t} \vert \{S^x_{t-1}\}] - \mathbb{E}[ \Phi_{t} \vert \{S^x_{t-1}\} ]}_{=\mathbf{0}} )]=0.
\end{eqnarray}
\end{lemma}

We next show Assumption \ref{ass:IVuncorrelatedWithYnoise} holds, which states that the instrumental variables are uncorrelated with the noise in the observations of the response variable. Below, $e_i$ denotes the column vector of all zeros except at the $i$th element which equals $1$.
\begin{lemma}\label{lem:ZYnoiseUncorrelated}
For every $i,j$, $\text{Cov}[Z_{ij}, Y^{''}_i] = 0$.
\newline \textbf{Proof:}
\begin{eqnarray}
\nonumber\text{Cov}[Z_{ij}, Y^{''}_i] &=& \mathbb{E}[Z_{ij} Y^{''}_i] - \mathbb{E}[Z_{ij}] \underbrace{\mathbb{E}[Y^{''}_i]}_{=0}\\
\nonumber &=& \mathbb{E}[\Phi_{t-1,ij} \{ \Pi_{t-1} (C_{t}-\overline{C}_{t})\}_i]\\
\nonumber &=& \mathbb{E}[\Phi_{t-1,ij} e^{\top}_i ~\Pi_{t-1} (C_{t}-\overline{C}_{t})]\\
\nonumber &=& \mathbb{E}\left[\mathbb{E}[ \Phi_{t-1,ij} e^{\top}_i~ \Pi_{t-1} (C_{t}-\overline{C}_{t}) \vert \{S^x_{t-1}\}]\right]\\
\nonumber &=& \mathbb{E}[\Phi_{t-1,ij} e^{\top}_i~ \Pi_{t-1} \underbrace{\mathbb{E}[C_{t}-\overline{C}_{t} \vert \{S^x_{t-1}\}]}_{=\vec{0}}]=0.
\end{eqnarray}
\end{lemma}

We next show Assumption \ref{ass:IVuncorrelatedWithXnoise} holds.
\begin{lemma}\label{lem:ZXnoiseUncorrelated}
For every $i$, $j$, and $l$, $\text{Cov}[Z_{ij}, X^{''}_{il}] = 0$.
\newline \textbf{Proof:}
\begin{eqnarray}
\nonumber\text{Cov}[Z_{ij}, X^{''}_{il}] &=& \mathbb{E}[Z_{ij} X^{''}_{il}]-\mathbb{E}[Z_{ij}] \underbrace{\mathbb{E}[ X^{''}_{il}]}_{=0}\\
\nonumber&=& \mathbb{E}[Z_{ij} (X^{'}_{il} - X_{il})]\\
\nonumber&=& \mathbb{E}[Z_{ij} e_i^{\top}(X^{'} - X)e_l] \label{eq:ZXerrorUncorrelated4}\\
\nonumber&=& \gamma \mathbb{E}[\Phi_{t-1,ij} e_i^{\top}\Pi_{t-1}(\mathbb{E}[\Phi_{t} \vert \{S^x_{t-1}\}] - \Phi_{t})e_l] \\
\nonumber&=& \gamma \mathbb{E}[\mathbb{E}[\Phi_{t-1,ij} e_i^{\top}\Pi_{t-1}(\mathbb{E}[\Phi_{t} \vert \{S^x_{t-1}\}] - \Phi_{t})e_l \vert \{S^x_{t-1}\}] ]\\
\nonumber&=& \gamma \mathbb{E}[\Phi_{t-1,ij} e_i^{\top}\Pi_{t-1}\mathbb{E}[\mathbb{E}[\Phi_{t} \vert \{S^x_{t-1}\}] - \Phi_{t} \vert \{S^x_{t-1}\}] e_l ]\\
\nonumber&=& \gamma \mathbb{E}[\Phi_{t-1,ij} e_i^{\top}\Pi_{t-1} \underbrace{(\mathbb{E}[\Phi_{t} \vert \{S^x_{t-1}\}] - \mathbb{E}[\Phi_{t} \vert \{S^x_{t-1}\}])}_{=0} e_l ]=0.
\end{eqnarray}
\end{lemma}
Finally, Assumption \ref{ass:IVsigmaFullRank} holds by Assumption \ref{ass:IVsigmaFullRankDP}, and assumptions \ref{ass:IVuncorrelatedWithYnoiseLimit}, \ref{ass:IVuncorrelatedWithXnoiseLimit}, and \ref{ass:IVsigmaFullRankLimit} follow trivially from assumptions \ref{ass:IVsigmaFullRankDP}, \ref{ass:OffPolicy}, \ref{ass:IVuncorrelatedWithYnoiseLimitLLN}, and \ref{ass:IVuncorrelatedWithXnoiseLimitLLN} by the law of large numbers (see Appendix \ref{app:IV}). Therefore Theorem \ref{thm:IVconsistent} follows. \qed
\end{theorem}

We proceed by showing that instrumental variables Bellman error minimization, least-squares projected Bellman error minimization, and instrumental variables projected Bellman error minimization methods, are equivalent.

\section{Equivalence of Instrumental Variables Bellman Error Minimization and Projected Bellman Error Minimization Methods}\label{sec:PolicyEvaluationAlgsEquivalent}
The following theorem points out the relationship among the three estimators in equations \eqref{eq:LSAPImatrix}, \eqref{eq:LSAPIprojected}, and \eqref{eq:IVProjectedBellmanErrorMinimization1}.

\begin{theorem}\label{thm:equivalencyOfPolicyEvaluationAlgs}
Under Assumption \ref{ass:PhiFullRank}, the following policy evaluation algorithms are equivalent:
\begin{itemize}
\item{Instrumental Variables Bellman Error Minimization (LSTD) in equation \eqref{eq:LSAPImatrix}:
\begin{eqnarray}
\nonumber&&\hat{\theta} = [(\Phi_{t-1})^{\top} (\Phi_{t-1} - \gamma {\Phi}_{t})]^{-1}(\Phi_{t-1})^{\top} C_{t}.
\end{eqnarray}}
\item{Least-Squares Projected Bellman Error Minimization (Least-Squares Fixed Point Approximation) in equation \eqref{eq:LSAPIprojected}:
\begin{eqnarray}
\nonumber&&\hat{\theta}= \left[\left(\Pi_{t-1} (\Phi_{t-1} - \gamma \Phi_{t})\right)^{\top} \left(\Pi_{t-1} (\Phi_{t-1} - \gamma \Phi_{t})\right)\right]^{-1} \left(\Pi_{t-1} (\Phi_{t-1} - \gamma \Phi_{t})\right)^{\top} \Pi_{t-1} C_{t}.
\end{eqnarray}}
\item{Instrumental Variables Projected Bellman Error Minimization in equation \eqref{eq:IVProjectedBellmanErrorMinimization1}:
\begin{eqnarray}
\nonumber&&\hat{\theta} = \left[ (\Phi_{t-1})^{\top} \Pi_{t-1} (\Phi_{t-1} - \gamma \Phi_{t}) \right]^{-1} (\Phi_{t-1})^{\top} \Pi_{t-1} C_{t}.
\end{eqnarray}}
\end{itemize}
\textbf{Proof:} We first show equations \eqref{eq:LSAPImatrix} and \eqref{eq:IVProjectedBellmanErrorMinimization1} are equivalent. Starting with equation \eqref{eq:IVProjectedBellmanErrorMinimization1} and recalling $\Pi_{t-1} = \Phi_{t-1} ((\Phi_{t-1})^{\top} \Phi_{t-1})^{-1} (\Phi_{t-1})^{\top}$, we can write
\begin{eqnarray}
&& \left[ (\Phi_{t-1})^{\top} \Pi_{t-1} (\Phi_{t-1} - \gamma \Phi_{t}) \right]^{-1} (\Phi_{t-1})^{\top} \Pi_{t-1} C_{t} \nonumber\\
&=& \left[ (\Phi_{t-1})^{\top}  (\underbrace{\Pi_{t-1} \Phi_{t-1}}_{\Phi_{t-1}} - \gamma \Pi_{t-1} \Phi_{t}) \right]^{-1} \underbrace{(\Phi_{t-1})^{\top} \Phi_{t-1} ((\Phi_{t-1})^{\top} \Phi_{t-1})^{-1}}_{I_{k \times k}} (\Phi_{t-1})^{\top} C_{t} \nonumber\\
&=& \left[ (\Phi_{t-1})^{\top}  \Phi_{t-1} - \gamma (\Phi_{t-1})^{\top}  \Pi_{t-1} \Phi_{t} \right]^{-1} (\Phi_{t-1})^{\top} C_{t} \nonumber\\
&=& \left[ (\Phi_{t-1})^{\top}  \Phi_{t-1} - \gamma \underbrace{(\Phi_{t-1})^{\top}  \Phi_{t-1} ((\Phi_{t-1})^{\top} \Phi_{t-1})^{-1}}_{I_{k \times k}} (\Phi_{t-1})^{\top} \Phi_{t} \right]^{-1} (\Phi_{t-1})^{\top} C_{t} \nonumber\\
&=& \left[ (\Phi_{t-1})^{\top}  \Phi_{t-1} - \gamma (\Phi_{t-1})^{\top} \Phi_{t} \right]^{-1} (\Phi_{t-1})^{\top} C_{t} \nonumber\\
&=& \left[ (\Phi_{t-1})^{\top} ( \Phi_{t-1} - \gamma \Phi_{t} ) \right]^{-1} (\Phi_{t-1})^{\top} C_{t}. \nonumber
\end{eqnarray}
Hence equations \eqref{eq:LSAPImatrix} and \eqref{eq:IVProjectedBellmanErrorMinimization1} are equivalent. Next, we show equations \eqref{eq:LSAPImatrix} and \eqref{eq:LSAPIprojected} are equivalent.
We start by writing
\begin{eqnarray}
&&(\Phi_{t-1} - \gamma \Phi_{t})^{\top} \Pi_{t-1} (\Phi_{t-1} - \gamma \Phi_{t})\nonumber\\
&=& (\Phi_{t-1} - \gamma \Phi_{t})^{\top} \Pi_{t-1} (\Phi_{t-1} - \gamma \Phi_{t})(\Phi_{t-1} - \gamma \Phi_{t})^{\top} \Phi_{t-1} [(\Phi_{t-1})^{\top} \Phi_{t-1}]^{-1} (\Phi_{t-1})^{\top} (\Phi_{t-1} - \gamma \Phi_{t})\nonumber\\
&=& (\Phi_{t-1} - \gamma \Phi_{t})^{\top} \Pi_{t-1} (\Phi_{t-1} - \gamma \Phi_{t})\nonumber\\
\implies &&(\Phi_{t-1} - \gamma \Phi_{t})^{\top} \Phi_{t-1} [(\Phi_{t-1})^{\top} \Phi_{t-1}]^{-1} \underbrace{(\Phi_{t-1})^{\top} (\Phi_{t-1} - \gamma \Phi_{t}) [(\Phi_{t-1})^{\top} (\Phi_{t-1} - \gamma \Phi_{t})]^{-1}}_{I_k} (\Phi_{t-1})^{\top} \nonumber\\
&=& (\Phi_{t-1} - \gamma \Phi_{t})^{\top} \Pi_{t-1} (\Phi_{t-1} - \gamma \Phi_{t})[(\Phi_{t-1})^{\top}(\Phi_{t-1} - \gamma \Phi_{t})]^{-1}(\Phi_{t-1})^{\top} \nonumber\\
\implies &&(\Phi_{t-1} - \gamma \Phi_{t})^{\top} \underbrace{\Phi_{t-1} [(\Phi_{t-1})^{\top} \Phi_{t-1}]^{-1}(\Phi_{t-1})^{\top}}_{\Pi_{t-1}} \nonumber\\
&=& (\Phi_{t-1} - \gamma \Phi_{t})^{\top} \Pi_{t-1} (\Phi_{t-1} - \gamma \Phi_{t})[(\Phi_{t-1})^{\top} (\Phi_{t-1} - \gamma \Phi_{t})]^{-1}(\Phi_{t-1})^{\top}(\Pi_{t-1} (\Phi_{t-1} - \gamma \Phi_{t}))^{\top} \nonumber\\
&=& (\Phi_{t-1} - \gamma \Phi_{t})^{\top} \Pi_{t-1} (\Phi_{t-1} - \gamma \Phi_{t})[(\Phi_{t-1})^{\top} (\Phi_{t-1} - \gamma \Phi_{t})]^{-1}(\Phi_{t-1})^{\top}\label{eq:IVtheoremProjectedBellman6}\\
\implies && [(\Phi_{t-1} - \gamma \Phi_{t})^{\top} \Pi_{t-1} (\Phi_{t-1} - \gamma \Phi_{t})]^{-1} (\Pi_{t-1} (\Phi_{t-1} - \gamma \Phi_{t}))^{\top} \nonumber\\
&=& \underbrace{[(\Phi_{t-1} - \gamma \Phi_{t})^{\top} \Pi_{t-1} (\Phi_{t-1} - \gamma \Phi_{t})]^{-1}(\Phi_{t-1} - \gamma \Phi_{t})^{\top} \Pi_{t-1} (\Phi_{t-1} - \gamma \Phi_{t})}_{I_k}\nonumber\\
\implies &&[(\Phi_{t-1})^{\top} (\Phi_{t-1} - \gamma \Phi_{t})]^{-1}(\Phi_{t-1})^{\top}\nonumber\\
\implies && [(\Phi_{t-1} - \gamma \Phi_{t})^{\top} \underbrace{\Pi_{t-1}}_{(\Pi_{t-1})^{\top} \Pi_{t-1}} (\Phi_{t-1} - \gamma \Phi_{t})]^{-1} (\Pi_{t-1} (\Phi_{t-1} - \gamma \Phi_{t}))^{\top} \nonumber\\
&=& [(\Phi_{t-1})^{\top} (\Phi_{t-1} - \gamma \Phi_{t})]^{-1}(\Phi_{t-1})^{\top}\nonumber\\
\implies && [(\Phi_{t-1} - \gamma \Phi_{t})^{\top} (\Pi_{t-1})^{\top} \Pi_{t-1} (\Phi_{t-1} - \gamma \Phi_{t})]^{-1} (\Phi_{t-1} - \gamma \Phi_{t})^{\top} \underbrace{(\Pi_{t-1})^{\top}}_{(\Pi_{t-1})^{\top} \Pi_{t-1}} \nonumber\\
&=& [(\Phi_{t-1})^{\top} (\Phi_{t-1} - \gamma \Phi_{t})]^{-1}(\Phi_{t-1})^{\top}\nonumber\\
\implies && [(\Phi_{t-1} - \gamma \Phi_{t})^{\top} (\Pi_{t-1})^{\top} \Pi_{t-1} (\Phi_{t-1} - \gamma \Phi_{t})]^{-1} (\Phi_{t-1} - \gamma \Phi_{t})^{\top} (\Pi_{t-1})^{\top} \Pi_{t-1} \nonumber\\
&=& [(\Phi_{t-1})^{\top} (\Phi_{t-1} - \gamma \Phi_{t})]^{-1}(\Phi_{t-1})^{\top}\nonumber\\
\implies && \left[\left(\Pi_{t-1} (\Phi_{t-1} - \gamma \Phi_{t})\right)^{\top} \left(\Pi_{t-1} (\Phi_{t-1} - \gamma \Phi_{t})\right)\right]^{-1} \left(\Pi_{t-1} (\Phi_{t-1} - \gamma \Phi_{t})\right)^{\top} \Pi_{t-1} \nonumber\\
&=& [(\Phi_{t-1})^{\top} (\Phi_{t-1} - \gamma \Phi_{t})]^{-1}(\Phi_{t-1})^{\top}\nonumber\\
\implies && \left[\left(\Pi_{t-1} (\Phi_{t-1} - \gamma \Phi_{t})\right)^{\top} \left(\Pi_{t-1} (\Phi_{t-1} - \gamma \Phi_{t})\right)\right]^{-1} \left(\Pi_{t-1} (\Phi_{t-1} - \gamma \Phi_{t})\right)^{\top} \Pi_{t-1} C_{t}\nonumber\\
&=& [(\Phi_{t-1})^{\top} (\Phi_{t-1} - \gamma \Phi_{t})]^{-1}(\Phi_{t-1})^{\top} C_{t}.\nonumber
\end{eqnarray}
Equation \eqref{eq:IVtheoremProjectedBellman6} uses the fact that $(\Pi_{t-1})^{\top}= \Pi_{t-1}$.  Hence equations \eqref{eq:LSAPImatrix}, \eqref{eq:LSAPIprojected}, and \eqref{eq:IVProjectedBellmanErrorMinimization1} are equivalent. \qed
\end{theorem}
The results in Section \ref{sec:NumericalWork} illustrate that the Least-squares Bellman error minimization cannot be equivalent to the Bellman error minimization with instrumental variables, and thus the two others.

\section{Direct Policy Search}\label{sec:DirectPolicySearch}
An alternative to Bellman error minimization for finding the regression vector $\theta$ is direct policy search.  
Consider policies parameterized by $\theta$ of the form in equation \eqref{X:Def}, in which the post-decision value function has been replaced by the linear model $\phi(S^x)^{\top} \theta$. In contrast to policy iteration or value iteration methods, the goal in direct policy search is not necessarily to find a value function which is close (with respect to some norm) to the true value function; the objective is to find a parameter vector $\theta$ for which the policy $X^\pi(s|\theta)$ performs well, i.e., it solves the following stochastic optimization problem
\begin{eqnarray}
\label{eq:DirectPolicySearchObjective}&\displaystyle\max_{\theta} F(\theta) = &\mathbb{E} \sum_{t=0}^T \gamma^t C(S_t,X^\pi(S_t|\theta)),
\end{eqnarray}
given the policy $X^\pi(S_t|\theta)$. Solving this problem becomes challenging, particularly as the dimension of $\theta$ grows. Furthermore, the optimization problem given by equation \eqref{eq:DirectPolicySearchObjective} is typically non-convex and non-separable. Classic stochastic optimization algorithms can be used to sequentially choose policies to simulate. When the dimension of $\theta$ is small, the Knowledge Gradient for Continuous Parameters (KGCP) policy has been shown to work well for efficiently optimizing $\theta$, e.g., see \citep{ScFrPo11}. This approach is explained in the following subsection.

Note that, in direct policy search, we only need to consider features which are a function of the decisions. In addition, $F(\theta)$ is nonconvex, lacks easily computable derivatives and can only be simulated (introducing noise), making the search for an optimal $\theta$ fairly difficult for higher dimensional problems.  In our experiments, $\theta$ is limited to two or three dimensions.



\subsection{The Knowledge Gradient for Direct Policy Search}
For a given value of $\theta^{i}$, we can obtain a noisy observation of the objective in equation \eqref{eq:DirectPolicySearchObjective} by simulating
\begin{eqnarray}
\nonumber&&\hat{V}^{\pi}(S_0)=C_0(S_0,X^\pi(S_0|\theta^{i})) + \gamma C_1(S_1,X^\pi(S_1|\theta^{i})) + \gamma^2 C_2(S_2,X^\pi(S_2|\theta^{i})) + \cdots.
\end{eqnarray}
The KGCP policy for optimizing $\theta$ treats the objective function $\mu(\theta)$ as a Gaussian process regression. It then combines a model of $\mu(\theta)\eqdef \hat{V}^{\pi}(S)$ with a criterion which chooses the next value of $\theta$ for which a noisy observation of $\mu(\theta)$ will be simulated. The KGCP quantifies how much the maximum of the objective is expected to increase by getting an additional noisy observation of $\mu(\theta)$ at a particular value of $\theta$.

More formally, let $\mathcal{F}^{n}$ be the sigma-algebra generated by $\theta^0,\cdots,\theta^{n-1}$ and the corresponding noisy observations of $\mu(\theta^0),\cdots,\mu(\theta^{n-1})$. 
Denote the updated expected values of $\mu$ at $\theta^{i}$, conditioned on $\mathcal{F}^{n}$, by ${\mu}^n(\theta^{i})$. 
The KGCP is defined as
\begin{eqnarray}
\label{eq:approxKG}&&\bar{\nu}^{KG,n}(\theta) \eqdef \mathbb{E}\left[\max_{i=0,..,n} {\mu}^{n+1}(\theta^i) \middle \vert \mathcal{F}^n, \theta^n = \theta\right] - \max_{i=0,..,n} {\mu}^{n}(\theta^i)|_{\theta^n=\theta} .
\end{eqnarray}
In the Gaussian process regression framework, ${\mu}^{n+1}(\theta)$ given $\mathcal{F}^n$ is normally distributed for each value of $\theta$, and the KGCP can be calculated exactly (see \citep{ScFrPo11} or its description in \cite{PoRy2012}[Chapter 16]). Thus, the next sampling decision is chosen to maximize the KGCP,
\begin{eqnarray}
\nonumber&&\theta^n \in \text{arg}\max_{\theta} \bar{\nu}^{KG,n}(\theta).
\end{eqnarray}
After $N$ observations, the implementation decision $\theta^n$ is chosen by maximizing the regression function,
\begin{equation}
\theta^{\star} \in \argmax_{\theta} {\mu}^{N}(\theta).\nonumber
\end{equation}
The value of $\theta$ which maximizes equation \eqref{eq:DirectPolicySearchObjective} produces the best policy within the class of polices, $X^\pi(S_t|\theta)$.  KGCP has been proven to converge asymptotically to the optimal value of $\theta$ (see \cite{ScFrPo11}).

\section{Benchmark application: energy storage}\label{sec:Models}
\begin{figure}[htb]
\begin{center}
\includegraphics[width=9cm, height=7.5cm]{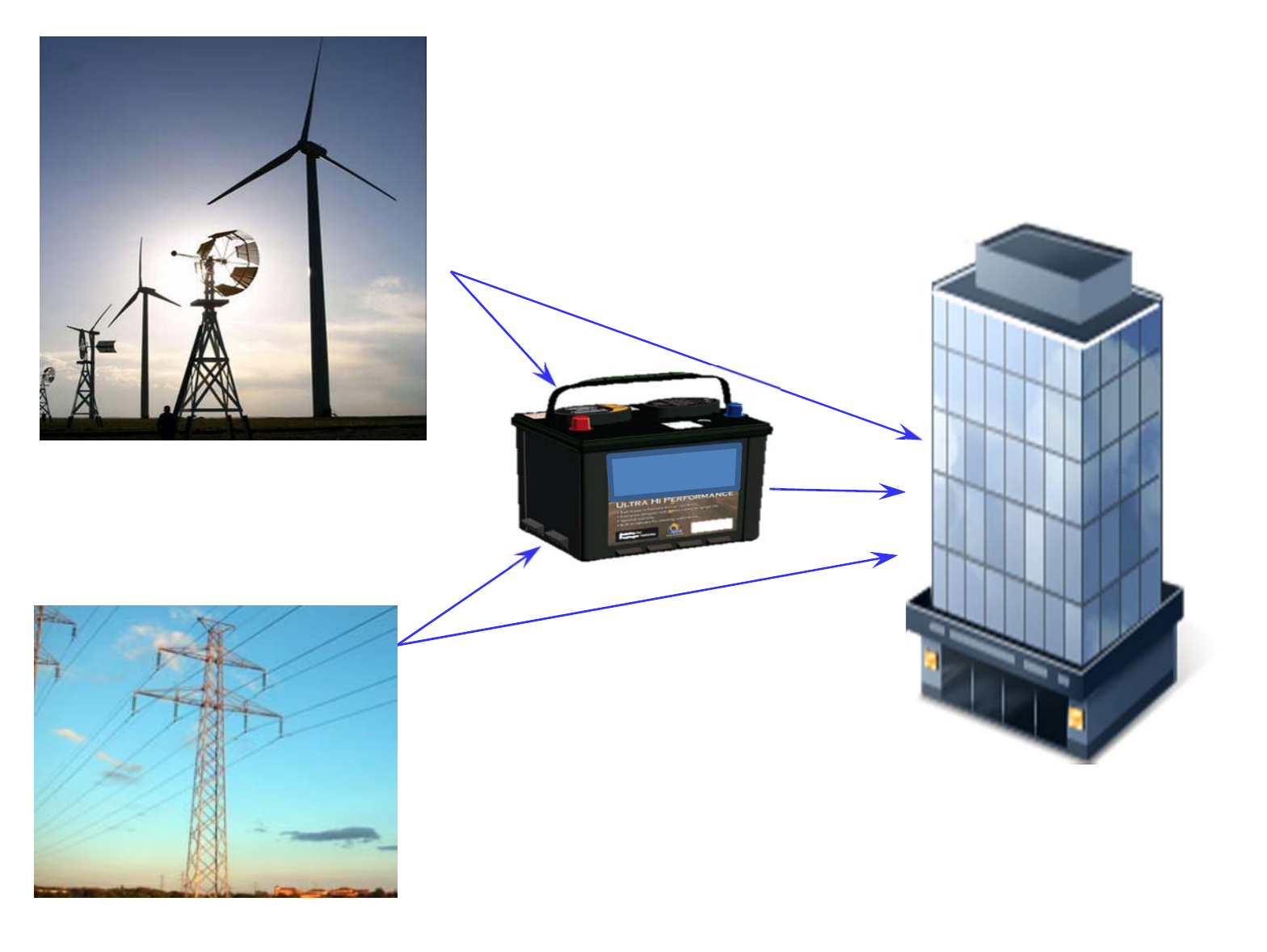}
\caption{The energy flow diagram.}\label{fig:energyFlows}
\end{center}
\end{figure}

Consider a power system as shown in Figure \ref{fig:energyFlows}, involving an intermittent energy supply, an electricity demand, a simple interconnecting grid, and a battery storage device. At time $t$, the energy flows in Figure \ref{fig:energyFlows} are given by the vector
\begin{eqnarray}
\nonumber&&x_{t}\eqdef \left(x^{WR}_t, x^{GR}_t, x^{RD}_t, x^{WD}_t, x^{GD}_t\right).
\end{eqnarray}
Here, $x_{t}^{IJ}$ denotes the amount of energy transferred from $I$ to $J$ at time step $t$. The superscript $W$ stands for energy source (wind), $D$ for demand, $R$ for storage, and $G$ for grid. All these entities are assumed to be nonnegative except $x^{GR}_{t}$. A negative value for $x^{GR}_{t}$ refers to selling the electricity from the storage to the grid.

Denote the total electricity demand (in $MWh$) over the time period starting at $t-\Delta t$ and ending at $t$, by $D_t$. At every time step, the  demand $D_t$ must be served through either the wind energy, available electricity from the storage device, or energy purchased from the grid,
\begin{eqnarray}
\label{eq:flowConstraint2}&& x^{GD}_{t}+\eta^{discharge} x^{RD}_{t}+x^{WD}_{t}=D_{t}.
\end{eqnarray}
Here, $0<\eta^{discharge}<1$ denotes the \emph{discharging efficiency rate}.

The wind energy generated during the time period $[t-\Delta t,~t)$, denoted by $E_{t}$, first flows to the electricity demand, i.e.,
\begin{eqnarray}
\label{eq:flowWind1}&&x^{WD}_{t}= \min\{E_t,~D_t\}.
\end{eqnarray}
The surplus is then sent to the storage device to be stored for the future use, i.e.,
\begin{eqnarray}
\label{eq:flowConstraint1} && x^{WR}_{t}+ x^{WD}_{t} = E_{t}.
\end{eqnarray}
This constraint assumes that the storage device can dissipate some part of the power received from the wind source upon receiving, if it is more than its capacity.

Let $R^{cap}$ indicate the total capacity of the storage. Define the constants $\Delta R^{min}$ and $\Delta R^{max}$ as the minimum and maximum fractions of the storage device that can be charged over $\Delta t$. For example, for a lead acid battery with a $C/10$ maximum charge and discharge rates, and $\Delta t = 15 \text{min}$, $\Delta R^{min}= -1/40$ and $\Delta R^{max}= 1/40$. To avoid charging or discharging the storage device faster than the permitted rates, $x^{GR}_t$ must satisfy:
\begin{eqnarray}
\label{eq:feasibleAction1}&&\frac{\Delta R^{min} R^{cap}}{\eta^{discharge}} \leq x^{GR}_t\leq \frac{\Delta R^{max} R^{cap}}{\eta^{charge}},
\end{eqnarray}
where $\eta^{charge}$ is the \emph{charging efficiency rate}. When selling the electricity from the storage to the grid is prohibited, the lower bound in equation \eqref{eq:feasibleAction1} must be replaced by $0$, i.e.,  $x^{GR}_t \geq 0$. Similarly, to ensure that the storage device is not discharged faster than allowed when sending energy from the storage device to the demand, we include:
\begin{eqnarray}
\label{eq:feasibleAction2}&&0 \leq x^{RD}_t \leq \Delta R^{max} R^{cap}.
\end{eqnarray}

By implementing the energy flow $X_{t}^{\pi}$, the next storage state can be computed as
\begin{eqnarray}
\label{Battery:R}&R_{t+\Delta t}(x_t) = \min\left\{R_t+\frac{\eta^{charge}( x^{GR}_t + x^{WR}_t) - x^{RD}_t}{R^{cap}}, 1\right\}.
\end{eqnarray}
In addition, we must have $R_{t+\Delta t}\geq R^{\min}$, where $R^{\min}$ is the minimum fraction of the capacity of the storage device that must remain full. For example, stationary lead-acid batteries with tubular plates, which are one of the lowest cost technologies for energy storage, should not be discharged below $20\%$ of their capacity, i.e., $R^{\min}=0.20$, e.g., see \citep{Br11}.

Constraints \eqref{eq:flowConstraint1}, \eqref{eq:feasibleAction1}, \eqref{eq:feasibleAction2}, \eqref{Battery:R}, in addition to standard measurability conditions, together define the set of admissible policies. Note that equations \eqref{eq:flowWind1}, \eqref{eq:flowConstraint1}, and \eqref{eq:flowConstraint2} effectively reduce the size of our action space from $5$ dimensions to $2$ dimensions.

The contribution function at every time step $t$ then is the dollar value of energy sold minus the amount bought from the grid,
\begin{eqnarray}
\label{Contributionfun}&&C(S_t,X_{t}^{\pi}) = P_{t} D_{t}-P_{t}\left(x^{GR}_{t} +x^{GD}_{t}\right).
\end{eqnarray}
The goal is to find a policy over the ergodic infinite horizon planning horizon which maximizes the accumulated expected discounted future rewards,
\begin{eqnarray}
\nonumber&& \max_\pi \mathbb{E}\left[\sum_{t=0}^{\infty}\gamma^t C\left(S_t, X^{\pi}_{t}\right) \right].
\end{eqnarray}
In the absence of the wind source and the electricity load, in which the storage device is used solely to buy/sell the electricity from/to
the grid, the model is referred to as the battery arbitrage problem.

Before seeking an optimal policy, it is necessary to define $E_{t},~D_{t}$, and prices $P_{t}$ over time steps. Uncertainties in the energy supply (wind), electricity demand, and prices are modeled through stochastic processes, as described in the subsequent subsections.

\subsection{Wind Power}\label{sec:Wind}
The power output from the wind turbine over $[t,~t+\Delta t)$, measured in $MWh$, is computed as below (see  \citep{Wind:Sustainable:MacKayBook} page $263$):
\begin{eqnarray}
\label{equation:wind:power}&&E_{t}=\frac{10^{-8}}{36}\times 0.5~ C_{p} ~\rho~ 50^{2}\pi~ W_{t}^{3}~\Delta t.
\end{eqnarray}
Here, $\rho$ is the density of air which is around $\rho=1.225~ kg/m^{3}$. The term $50^{2}\pi~m^{2}$ approximates the area in square meters swept by the rotor blades of the turbine. The power coefficient is around $C_p=0.45$. In equation \eqref{equation:wind:power}, $W_{t}$ denotes the wind speed measured in meters per second and $\Delta t$ is stated in seconds. Velocity of the wind, $W_{t}$, in $m/s$, equals
\begin{eqnarray}
\nonumber&&W_{t}=\left(Y_{t}+\mathbb{E}\left[\sqrt{W_t}\right]\right)^{2}.
\end{eqnarray}
where $Y_t$ is the de-meaned square root of the wind speeds. Similar to \cite{BrKaMu84}, we let $Y_t$ evolve by an AR(1) model:
\begin{eqnarray}
\nonumber&&Y_t = \phi_1 Y_{t-\Delta t} +  \epsilon _t,~~~~~\epsilon_t \sim \mathcal{N}(0, \sigma_{\epsilon}).
\end{eqnarray}
Using $15$-min data from the wind speeds at Maryneal, Texas and applying the Yule-Walker equations, e.g., see \citep{Ca04} to fit the above model, we obtained $\mathbb{E}[\sqrt{W_t}] = 1.4781$, $\phi_1 = 0.7633$, $\sigma_{\epsilon} = 0.4020$. Note that other models have also been used in the literature; for example, \cite{ChPeBaCh10} suggest using a more general ARIMA model for the square root of the wind speed.

\subsection{Electricity Prices}\label{sec:ElectricityPrices}
\begin{figure}[htb]
\begin{center}
\includegraphics[width=9cm, height=7.5cm]{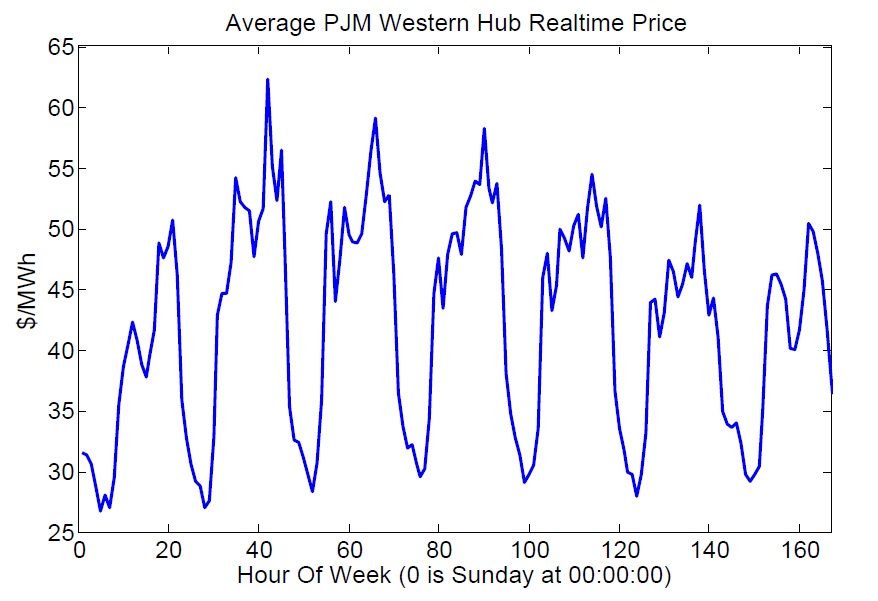}
\caption{The average PJM Real-Time price at the Western Hub as a function of the hour of the week. }\label{fig:PJMpricesRealTimeMeanHrOfWeek}
\end{center}
\end{figure}

Figure \ref{fig:PJMpricesRealTimeMeanHrOfWeek} illustrates the electricity spot prices at the PJM Western Hub for one week. While the average of these real-time prices over 2009-2010 is $\$42.11$ per $MWh$, the prices are occasionally negative and have a maximum of $\$362.90$ per $MWh$. Figure \ref{fig:PJMpricesRealTimeMeanHrOfWeek} shows that the prices are lowest at night; they begin to increase around $5~am$ and are typically the highest in the evening around $6~pm$.

In addition to seasonality, mean reversion and heavy-tailed behaviour have been considered as important properties of energy prices, (see  \cite{EyWo03}). Therefore, as in \cite{CaFi05}, here we adopt a mean reverting jump diffusion model with seasonality for real-time electricity prices in $\$/MWh$,
\begin{eqnarray}
\label{eq:logPrices}&& P_{t}=\exp\left(Y^{ds}_{t}+Y^{s}_{t}\right)-c,
\end{eqnarray}
where the constant parameter $c$ equals one minus the minimum value of $P_t$ in the data set, i.e., $c=1-\min_{t}~P_{t}$. In equation \eqref{eq:logPrices}, the seasonal component $Y^{s}_t$ is a deterministic periodic function of $t$, including an hour-of-week and a
month-of-year components. The deseasonalized log prices are modeled using the following discretized mean reverting jump diffusion model:
\begin{eqnarray}
\label{eq:jumpDiffusionDiscretized}&&Y^{ds}_t - Y^{ds}_{t-\Delta t} = \lambda(\mu - Y^{ds}_{t- \Delta t}) \Delta t + \sigma \sqrt{\Delta t}~ \epsilon_t + J_t,
\end{eqnarray}
where $\mu$ is the long term equilibrium price, $\lambda$ is the mean reversion rate, $\{\epsilon_{t+n \Delta t}\}_{n=0}^{N}$ are i.i.d. standard normal random variables, and $J_t$ is the jump over the interval $[t- \Delta t,~t)$.
As in \cite{CaFi05}, we model the jumps by the i.i.d. process,
\begin{eqnarray}
\nonumber&&J_t = \epsilon^{jump}_{t} \mathbf{1}\left(U_t < p^{jump}\right),
\end{eqnarray}
where $\epsilon^{jump}_t\sim N(0,\sigma^{jump})$ is the size of a jump, $U_t \sim \text{unif}(0,1)$, and $p^{jump}$ is the probability of a jump over a time interval of length $\Delta t$. In addition, jump sizes $\{\epsilon^{jump}_{t + n \Delta t}\}_{n=0}^{N}$ are i.i.d. normal random variables.

We identify the nonzero jumps (as in \cite{CaFi05}) by counting the number of times where the absolute value of the return is more than three times its standard deviation. The probability of a jump $p^{jump}$ equals then half of the fraction of times that jumps occur; we divide this estimate by two because most jumps are immediately followed by jumps in the opposite direction. The standard deviation of the magnitude of these jumps is considered as $\sigma^{jump}$. After estimating the jump terms, the parameters $\lambda$, $\mu$, and $\sigma$ can be estimated using least-squares linear regression on equation \eqref{eq:jumpDiffusionDiscretized}. For the electricity spot prices at the PJM Western Hub data set, we get $\sigma^{jump}=0.4229$, $p^{jump}=0.0170$, $\lambda=1800.9$, $\mu=4.1995$, $\sigma=11.0971$, and $c=27.2531$.



\subsection{Electricity Demand}\label{sec:ElectricityDemand}
\begin{figure}[htb]
\begin{center}
\includegraphics[width=7.5cm, height=6cm]{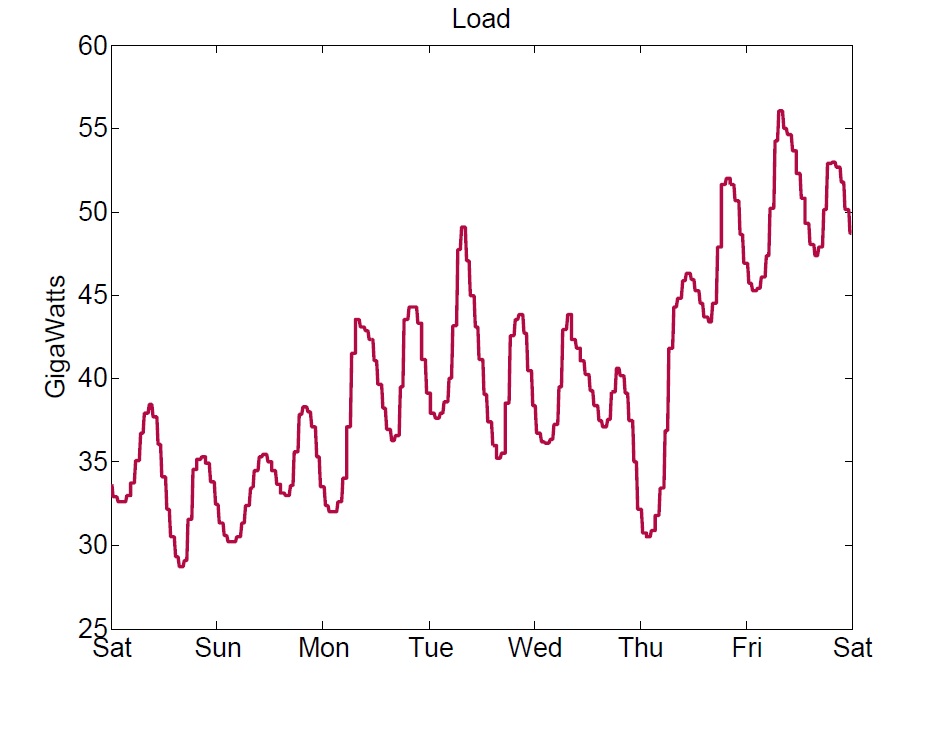}
\includegraphics[width=7.5cm, height=6cm]{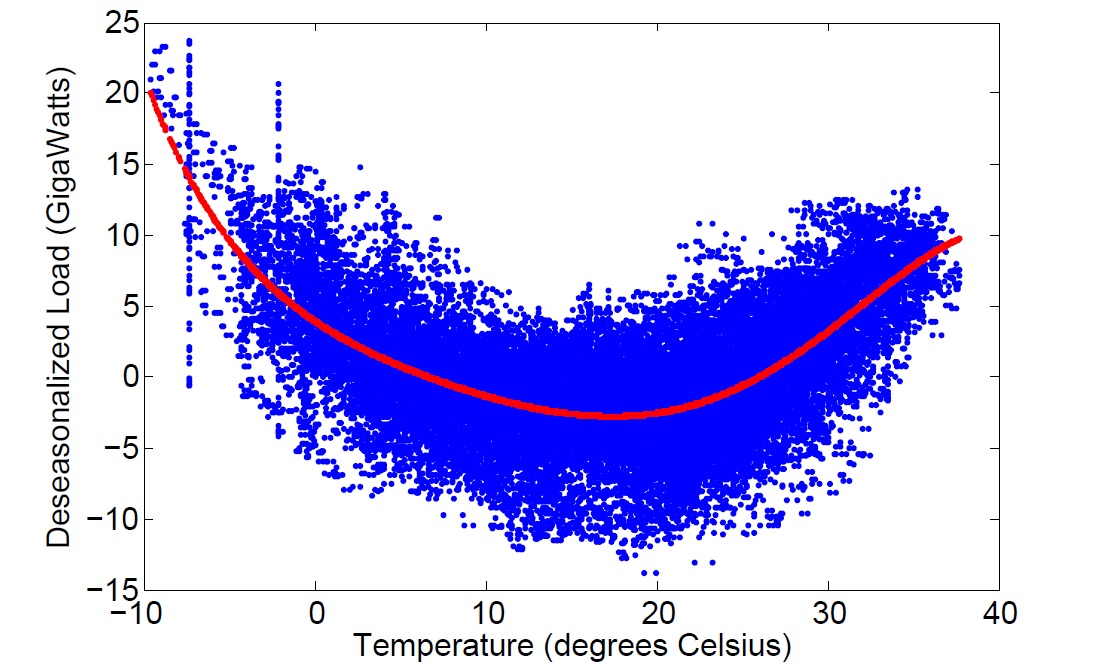}
\end{center}
\caption{Left plot: The total ERCOT load over the first full week of 2010. Right plot: Deseasonalized total ERCOT load versus temperature at Maryneal, TX over 2010.  There is one point for each 15 minute period of 2010.  The fitted fifth order polynomial regression is shown as well.}
\label{fig:loadDSversusTemp}
\end{figure}

\cite{EyWo03} outline typical models for residential, commercial, and industrial power demands. While industrial power demand is relatively stable, residential power demand is highly dependent upon the temperature and exhibits seasonal variations. These hourly and daily features can be clearly observed from the left plot in Figure \ref{fig:loadDSversusTemp}; the hourly ERCOT energy load starts ramping up in the morning around $5$~am and peaks in the evening around $6$~pm, although the patterns vary greatly based on the day of the week and the month of the year. Many approaches to forecast loads, such as an end-use model that incorporates appliances and customers, various regression models (based on temperature, time, and other factors), and heuristics made by experts, are summarized in \cite{FeGe05}. Similar to \cite{PiJe01} (and for the sake of simplicity) we adopt a Brownian motion with seasonality components for modeling the electricity demands in $MWh$,
\begin{eqnarray}
\label{eq:loadSeasonal}&&D_t = m^{hour}_{t} + m^{month}_{t} + D^{ds}_{t}.
\end{eqnarray}
Here, $m^{hour}_t$ and $m^{month}_t$ indicate the hour-of-week seasonal and the month-of-year seasonal components, and $D^{ds}_t$ is the deseasonalized load. The hour-of-week seasonal component $m^{hour}_t$ is calculated as the average load over each of the hours of the week, and
$m^{month}_t$ is computed as the average load over each of the months of a year. The deseasonalized load $D^{ds}_t$ evolves with
temperature (in degrees Celsius) $T_{t}$ at time $t$, as
\begin{eqnarray}
\label{eq:loadDeseasonalizedRegression}&&D^{ds}_{t} = \alpha_0 + \alpha_1 T_t + \alpha_2 T_t^2 + \alpha_3 T_t^3 + \alpha_4 T_t^4 + \alpha_5 T_t^5 + \epsilon_t,~~~~~\epsilon_t \sim \mathcal{N}(0, \sigma_{\epsilon})
\end{eqnarray}
where the coefficients $\alpha_{0},\alpha_{1},\cdots,\alpha_{5}$ can be estimated through least squares polynomial regression. The right plot in Figure \ref{fig:loadDSversusTemp} illustrates the estimated deseasonalized load $D^{ds}_{t}$ against temperature. For the hourly ERCOT energy load data we obtained, $\phi_1=0.9636$, $\sigma_{\epsilon}^{2} = 914870$.

\subsection{Stochastic Dynamic Programming Formulation}\label{sec:DynamicProgram}
We formulate this multi-stage stochastic programming problem through stochastic dynamic programming. The state variable, $S_t$, and the
post-decision state variable, $S^{x}_{t}$, are defined as,
\begin{eqnarray}
\nonumber&&S_t = \left(R_t,~E_t,~D_t,~P_t\right),\\
\nonumber&&S^{x}_{t} = \left(R_{t+\Delta t}(x),~ E_t,~ D_t,~ P_t\right).
\end{eqnarray}
Recall that $R_t$ is the fraction of the storage device that is full, $E_t$ is the current amount of wind energy, $D_t$ indicates the current energy demand, and $P_{t}$ is the current spot price of electricity selling to (or purchasing from) the electrical grid. 
The \emph{exogenous information process} is defined as the random changes in the state of the system, $W_{t+\Delta t}=\{\hat{E}_{t+\Delta t},\hat{D}_{t+\Delta t},\hat{P}_{t+\Delta t}\}$, which refer to exogenous changes in $E_t$, $D_t$ and $P_t$. These exogenous changes may be state dependent as well as time dependent.

The goal is then to find the policy, $X^\pi(S_t)$, which maximizes the expected discounted future rewards,
\begin{eqnarray}
\label{eq:DynamicProgrammingObjective} &\displaystyle\max_{\pi \in \Pi}& \mathbb{E}\left[\sum_{t=0}^{\infty}\gamma^{t}~ C\left(S_t, X^\pi(S_t)\right) \right].
\end{eqnarray}
The policy $X^\pi(S_t)$ may be stationary or time-dependent, which may be denoted by $X^\pi_t(S_t)$. Time-dependent policies can be explained by modifying the state variable through including time.

In the following section, we use the LSAPI algorithm, described in Section \ref{sec:APIA} with policy evaluation methods discussed in Section \ref{Sec:LSTD:PE}, as well as the direct policy search described in Section \ref{sec:DirectPolicySearch}, for finding policies, which are then compared to optimal policies so we can evaluate their performance.

\section{Numerical Experiments}\label{sec:NumericalWork}
Our main objective in this section is to assess the performance of the API variants and direct policy search. In subsection \ref{sec:Benchmark}, these algorithms are compared with $20$ discretized benchmark problems for energy storage operation management. The size of the benchmark problems considered are realistic relative to real world instances. Exact optimal stochastic dynamic programming solutions of these benchmark problems have been computed using the exact value iteration method.  Finally, we run approximate policy iteration on a problem with a state consisting of five continuous dimensions and actions consisting of five continuous dimensions.

Throughout, the time step is fixed to $\Delta t=15~\min$. We use a discount factor, $\gamma = 99.90\%$. We found that discount factors of $\gamma=99\%$ or smaller produce policies that are relatively myopic, and do not allow us to hold energy in storage for extended periods.

\subsection{Benchmark Problems with Discrete State Spaces}\label{sec:Benchmark}
We first consider finite, discretized state and action spaces with a fixed probability transition matrix. An exact solution for an infinite horizon problem can be found using the value iteration method, e.g., see \cite{Pu94}. We note that computing these optimal policies typically required approximately two {\it weeks} of CPU time, due primarily to the very high discount factor. In this method, $V^0(s)$ is initialized to a constant for all states $s$ in the state space, and at each iteration $n$, the algorithm updates the value function at each state using
\begin{eqnarray}
\nonumber&&V^n(s) = \text{max}_{x} \left\{ C(s, x) + \gamma \sum_{s'\in\mathcal S} V^{n-1}(s')\Pr(s'|s,x) \right\} , \quad \text{for every  } s \in \mathcal{S}.
\end{eqnarray}
The algorithm will converge to within $\epsilon$ of the true value functions (see \cite{Pu94}).

To establish a set of benchmark problems, we considered $20$ instances of the energy storage operation management problem in Section \ref{sec:Models}. These problems are summarized in Table \ref{tab:benchmark}, where problem type ``Full'' refers to the problem in Figure \ref{fig:energyFlows} with energy from wind source and the grid serves an electricity demand. The problem type ``BA'' refers to a battery arbitrage problem, which considers only trading between the storage and the grid to take advantage of price variations. We discretized the state space in the benchmark test problems and then created fixed probability transition matrices for the exogenous information process in order to create a true discrete process.
Table \ref{tab:benchmark} also reports how finely each state variable is discretized (the size of the state space for a particular problem is the product of each of the discretization levels). We then list the average maximum wind capacity divided by the load, the storage capacity divided by the load over an hour, the round trip efficiency (RTE) of the storage device, and the max charge and discharge rate of the storage device. For example, $C/10$ indicates that the storage device can be fully charged or discharged within $10$ hours. The transition matrix of the electricity prices was fitted using the PJM Western Hub real time prices (with and without time of day). The transition matrix of the load was fitted using the load of the PJM Mid-Atlantic Region (with time of day). The transition matrix for the wind was fit using data from wind speeds near the Sweetwater Wind Farm.  For Problems $1-16$ the state space is resource level, wind energy, and electricity price, i.e., $S_t = (R_t, E_t, P_t)$. For this experiment, time and demand are held fixed in order to keep the benchmark problems computationally tractable, as the exact value iteration, even for these simplified problems, requires approximately two weeks on a 3Ghz processor. For Problems $17-20$, the state space is the time-of-day, $\tau_t$, (1-96 corresponding to fifteen minute intervals in a day), the resource level, and the electricity price. The full state variable is then given by $S_t = (\tau_t, R_t, P_t)$.

\begin{table}[htb]
\begin{center}
\small{
\begin{tabular}{|c|c|c|c|c|c|c|c|c|c|c|}
\hline
               \multicolumn{2}{|c|}{Problem}             & \multicolumn{5}{c|}{Number of Discretization Levels}& \multicolumn{4}{c|}{Parameters} \\
\hline
Number & Type & Time & Resource & Price & Demand & Wind & Wind  & Storage & RTE & Charge Rate \\
       &      & $\tau_{t}$ & $R_{t}$ & $P_{t}$ & $D_{t}$ & $E_{t}$ & $\frac{E_{t}}{D_{t}}$ & $\frac{R_{t}}{D_{t}}$  & &  \\
\hline
1 & Full & 1 & 33 & 20 & 1 & 10  & 0.1 & 2.5 & .81 & C/10 \\
\hline
2 & Full & 1 & 33 & 20 & 1 & 10  & 0.1 & 2.5 & .81 & C/1 \\
\hline
3 & Full & 1 & 33 & 20 & 1 & 10  & 0.1 & 2.5 & .70 & C/10 \\
\hline
4 & Full & 1 & 33 & 20 & 1 & 10  & 0.1 & 2.5 & .70 & C/1 \\
\hline
5 & Full & 1 & 33 & 20 & 1 & 10  & 0.2 & 2.5 & .81 & C/10 \\
\hline
6 & Full & 1 & 33 & 20 & 1 & 10  & 0.2 & 2.5 & .81 & C/1 \\
\hline
7 & Full & 1 & 33 & 20 & 1 & 10  & 0.2 & 2.5 & .70 & C/10 \\
\hline
8 & Full & 1 & 33 & 20 & 1 & 10  & 0.2 & 2.5 & .70 & C/1 \\
\hline
9 & Full & 1 & 33 & 20 & 1 & 10  & 0.1 & 5.0 & .81 & C/10 \\
\hline
10 & Full & 1 & 33 & 20 & 1 & 10 & 0.1 & 5.0 & .81 & C/1 \\
\hline
11 & Full & 1 & 33 & 20 & 1 & 10 & 0.1 & 5.0 & .70 & C/10 \\
\hline
12 & Full & 1 & 33 & 20 & 1 & 10 & 0.1 & 5.0 & .70 & C/1 \\
\hline
13 & Full & 1 & 33 & 20 & 1 & 10 & 0.2 & 5.0 & .81 & C/10 \\
\hline
14 & Full & 1 & 33 & 20 & 1 & 10 & 0.2 & 5.0 & .81 & C/1 \\
\hline
15 & Full & 1 & 33 & 20 & 1 & 10 & 0.2 & 5.0 & .70 & C/10 \\
\hline
16 & Full & 1 & 33 & 20 & 1 & 1  & 0.2 & 5.0 & .70 & C/1 \\
\hline
17 & BA & 96 & 33 & 20 & 1 & 1 & - & - & .81 & C/10 \\
\hline
18 & BA & 96 & 33 & 20 & 1 & 1 & - & - & .81 & C/1 \\
\hline
19 & BA & 96 & 33 & 20 & 1 & 1 & - & - & .70 & C/10 \\
\hline
20 & BA & 96 & 33 & 20 & 1 & 1 & - & - & .70 & C/1 \\
\hline
\end{tabular}
}
\end{center}
\caption{Benchmark problems: the number of discretization levels for time (1=steady state), resource, price, load (1=deterministic) and wind.}
\label{tab:benchmark}
\end{table}

\begin{figure}[htb]
\begin{center}
\includegraphics[width=7.5cm, height=6cm]{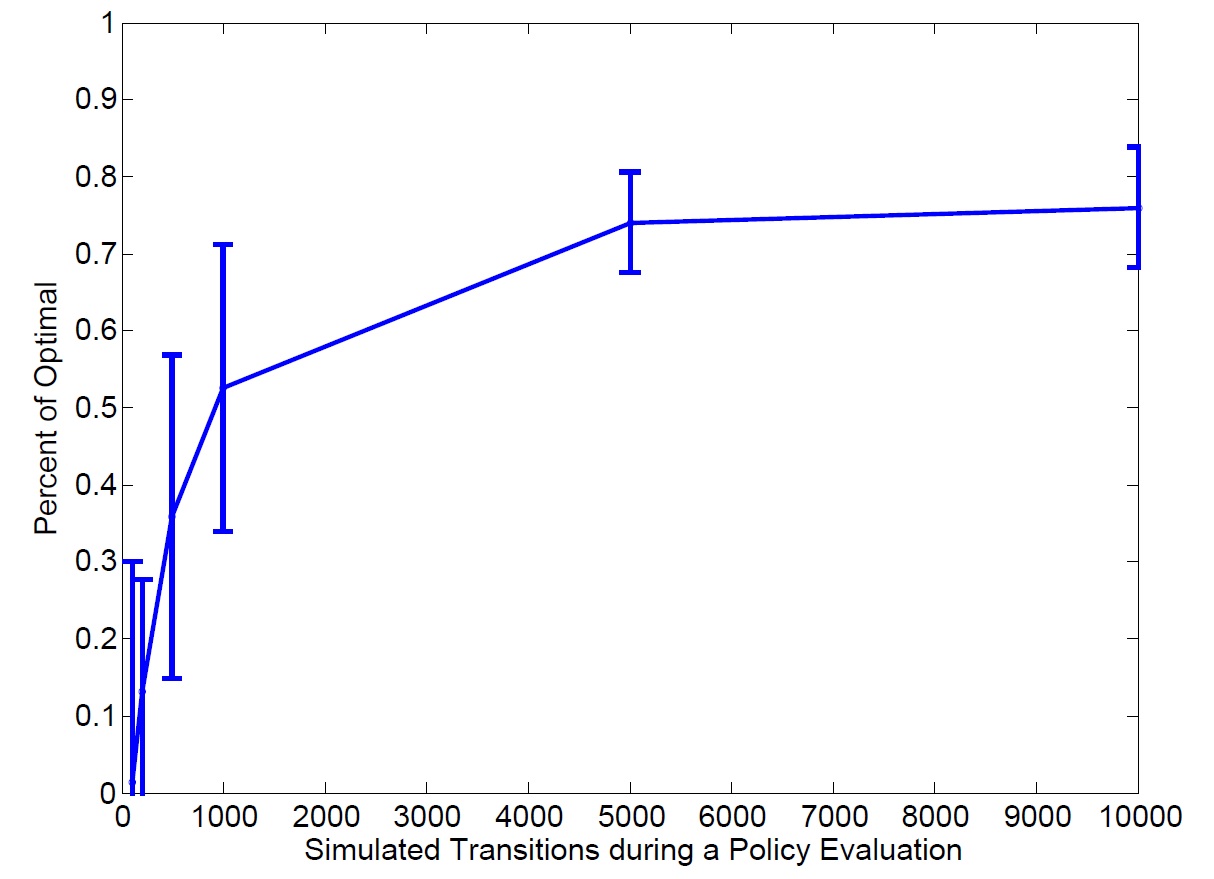}
\includegraphics[width=7.5cm, height=6cm]{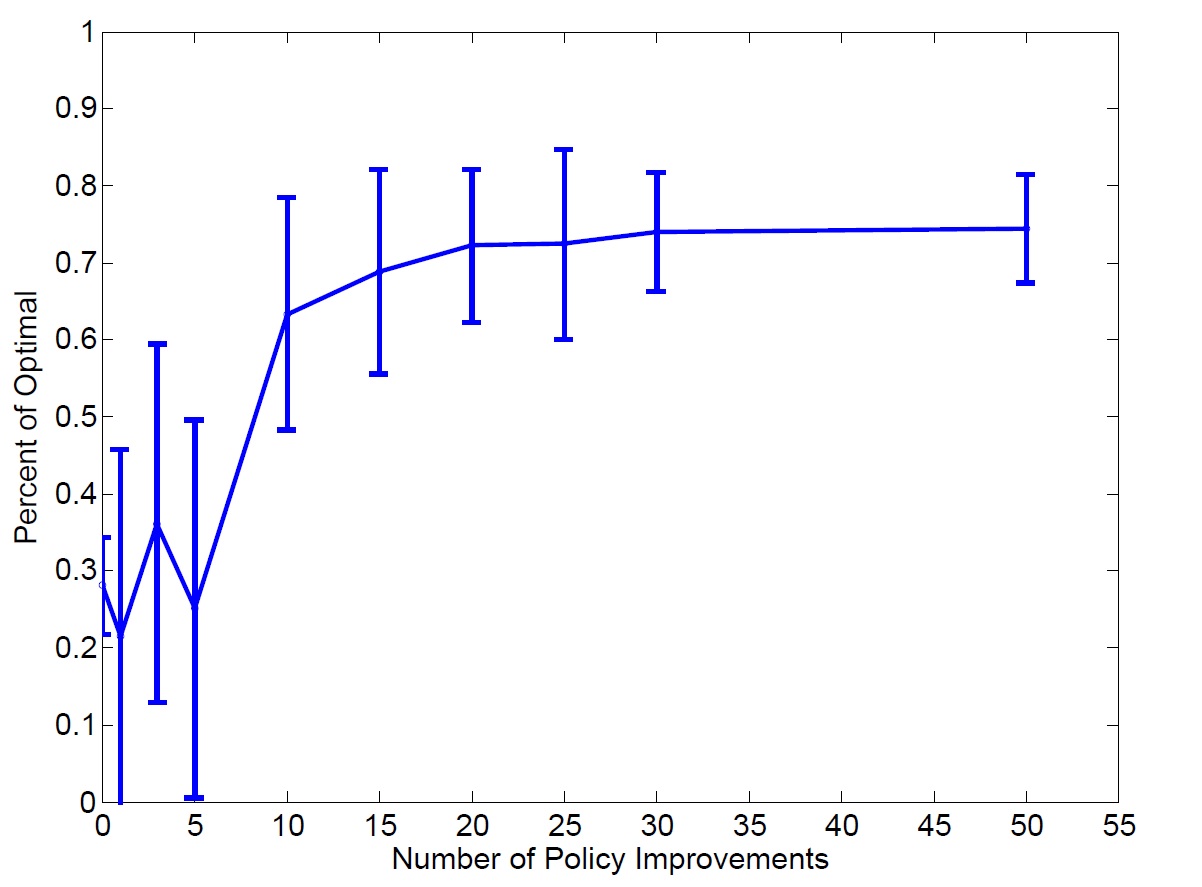}
\caption{Left Plot: Progress of the LSAPI method using instrumental variable Bellman error minimization for benchmark problem $17$ as $N$ varies. Right Plot: Progress of the LSAPI method using instrumental variable Bellman error minimization for benchmark problem $17$ as $M$ varies.}\label{fig:varyLengthPolicyEvaluationProblem17}
\end{center}
\end{figure}

To specify reasonable values for the maximum number of policy evaluation and policy improvement iterations, $M$ and $N$ in Figure \ref{fig:LSAPI}, we implemented the LSAPI method using instrumental variables Bellman error minimization on the $17^{th}$ benchmark problem several times. For this test problem, as illustrated in Figure \ref{fig:varyLengthPolicyEvaluationProblem17}, most of the improvement has occurred before $N=5000$ iterations of policy evaluations and $M=30$ iterations of policy improvement step. Thus, in the rest of this section, we fix $N=5000$ and $M=30$.

\subsection{API and Direct Policy Search Policies for the Benchmark Problems}\label{sec:BenchmarkComparisson}
\begin{figure}[htb]
\begin{center}
\includegraphics[width=6in,viewport=50 200 600 600]{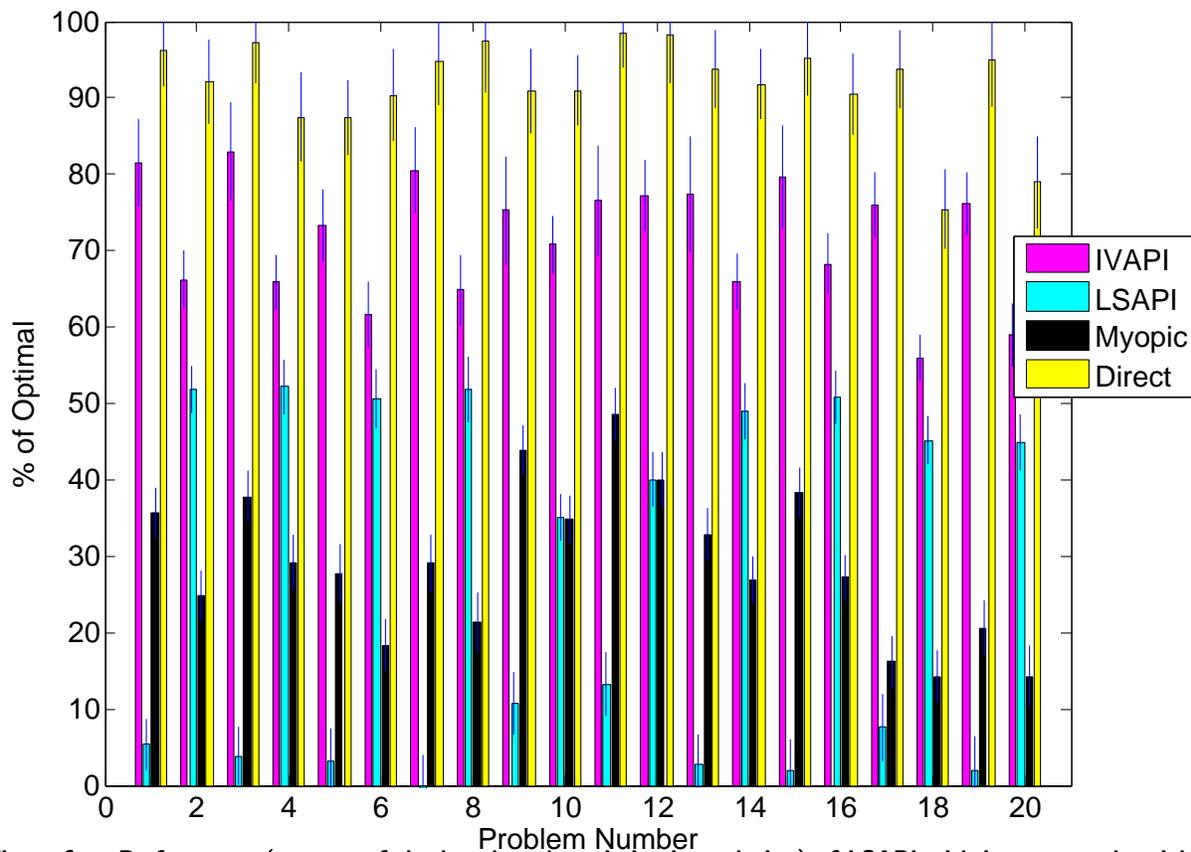}

\caption{Performance (percent of the benchmark optimization solution) of LSAPI with instrumental variables, LSAPI, the myopic policy and direct policy search.}\label{fig:perfAllProblemsOrder1and2}

\end{center}
\end{figure}

Figure \ref{fig:perfAllProblemsOrder1and2} compares the LSAPI method with instrumental variables Bellman error minimization (IVAPI), API with least-squares Bellman error minimization (LSAPI), the myopic policy (Myopic), and direct policy search based on KGCP (Direct). In this figure, quadratic basis functions are used to implement the IVAPI and LSAPI methods. The myopic policy discharges the storage device as quickly as possible and keeps it at its minimum level since then. The value of the myopic policy is still positive due to the wind source as well as the initial storage level.

We run each algorithm $100$ times. For each run of the algorithms, the final policies produced by each algorithm are then evaluated on the same sample paths, $\omega \in \Omega$, where $\omega$ is generated from the discretized exogenous information process.  We then record the average percent of optimal and the standard deviation of the average percent of optimal across the $100$ runs.  The average percentage of optimal for a policy $\pi$ is computed as
\begin{equation}
\text{\% of optimal} = \frac{1}{|\Omega|}\sum_{\omega \in \Omega} \frac{\hat{F}^{\pi}\left(S_{0}(\omega)\right)}{V^{\star}\left(S_0(\omega)\right)}, \nonumber
\end{equation}
where $\omega$ is a sample path of the randomness in the state transitions, and $S_0(\omega)$ is the starting state which has been randomly generated from a uniform distribution.  $\hat{F}^{\pi}(S_{0}(\omega))$ is a realization of value of the policy $\pi$ run on the sample path $\omega$, starting at the state $S_0(\omega)$, and $V^{\star}(S_0(\omega))$ is the true value of the optimal policy for state $S_0(\omega)$ computed using the exact value iteration method.

Similarly, we implemented direct policy search using KGCP $100$ times, and computed the average percent of optimal and its standard deviation.
To implement direct policy search using KGCP, we budget ourselves to simulating $50$ sequentially chosen policies, after which the KGCP algorithm must choose what it believes to be the best policy. 

Figure \ref{fig:perfAllProblemsOrder1and2} illustrates that IVAPI significantly outperforms LSAPI, but still underperforms the optimal policy by a wide margin. Direct policy search produced solutions that are on average $91.80$ percent of optimal, and are always at least $70\%$ of optimal for problems $1$ through  $16$. This suggests that for the benchmark problems on the application of interest, direct policy search is much more robust relative to the LSAPI methods. However, direct policy search quickly becomes intractable as the number of basis functions increases. In addition, choosing the search domain for direct policy search is another significant complication as the number of basis functions increases. We suggest using LSAPI to find good values of the regression parameters, and then apply direct policy search to improve the policy in the region of the fitted regression parameters. 

\begin{figure}[htb]
\begin{center}
    \subfigure[]{\label{fig:Problem1IVsamplePath3}\includegraphics[width=3in,viewport=50 200 550 550]{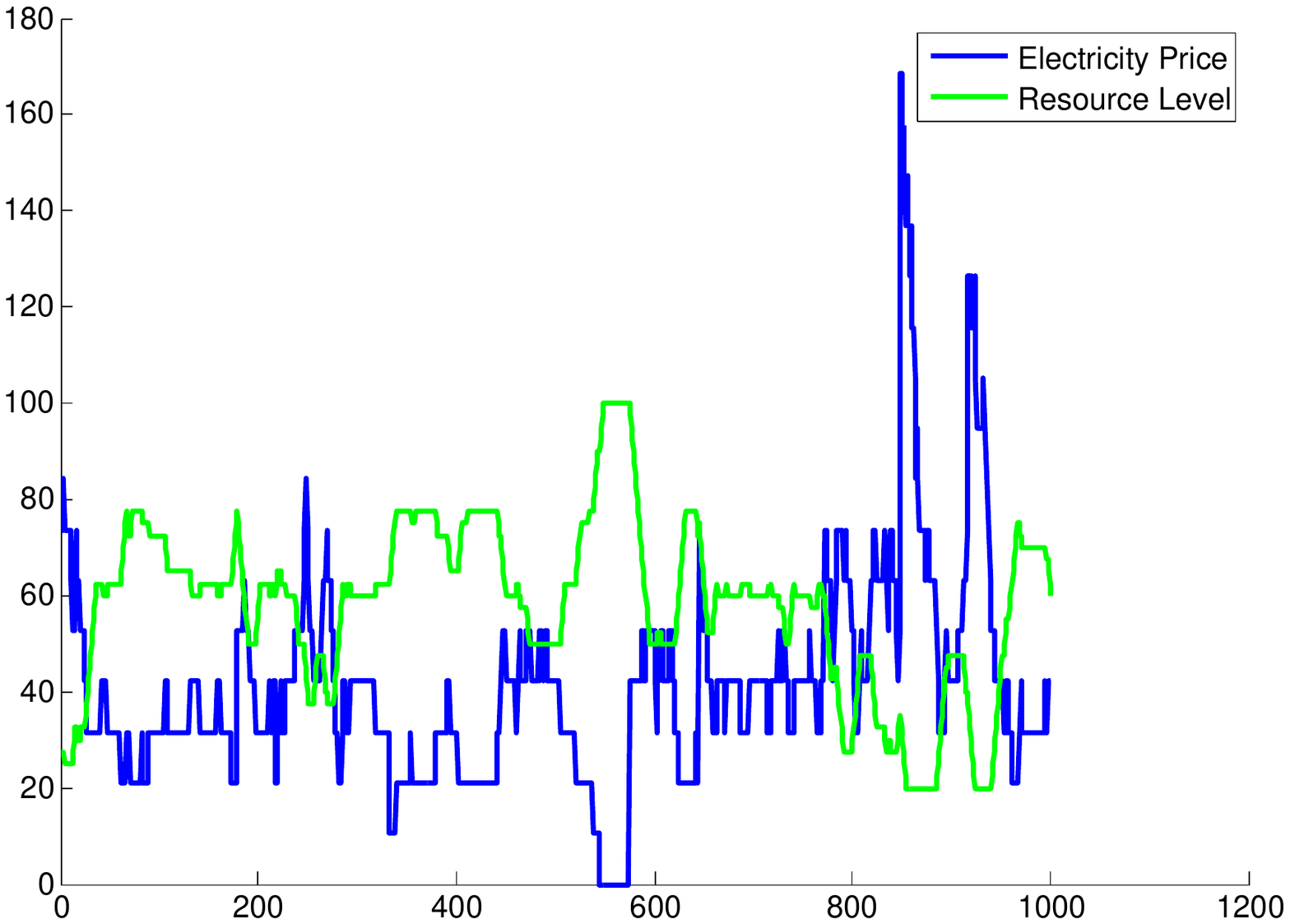}}
    \subfigure[] {\label{fig:Problem1IVhistogram}\includegraphics[width=3in,viewport=50 200 550 550]{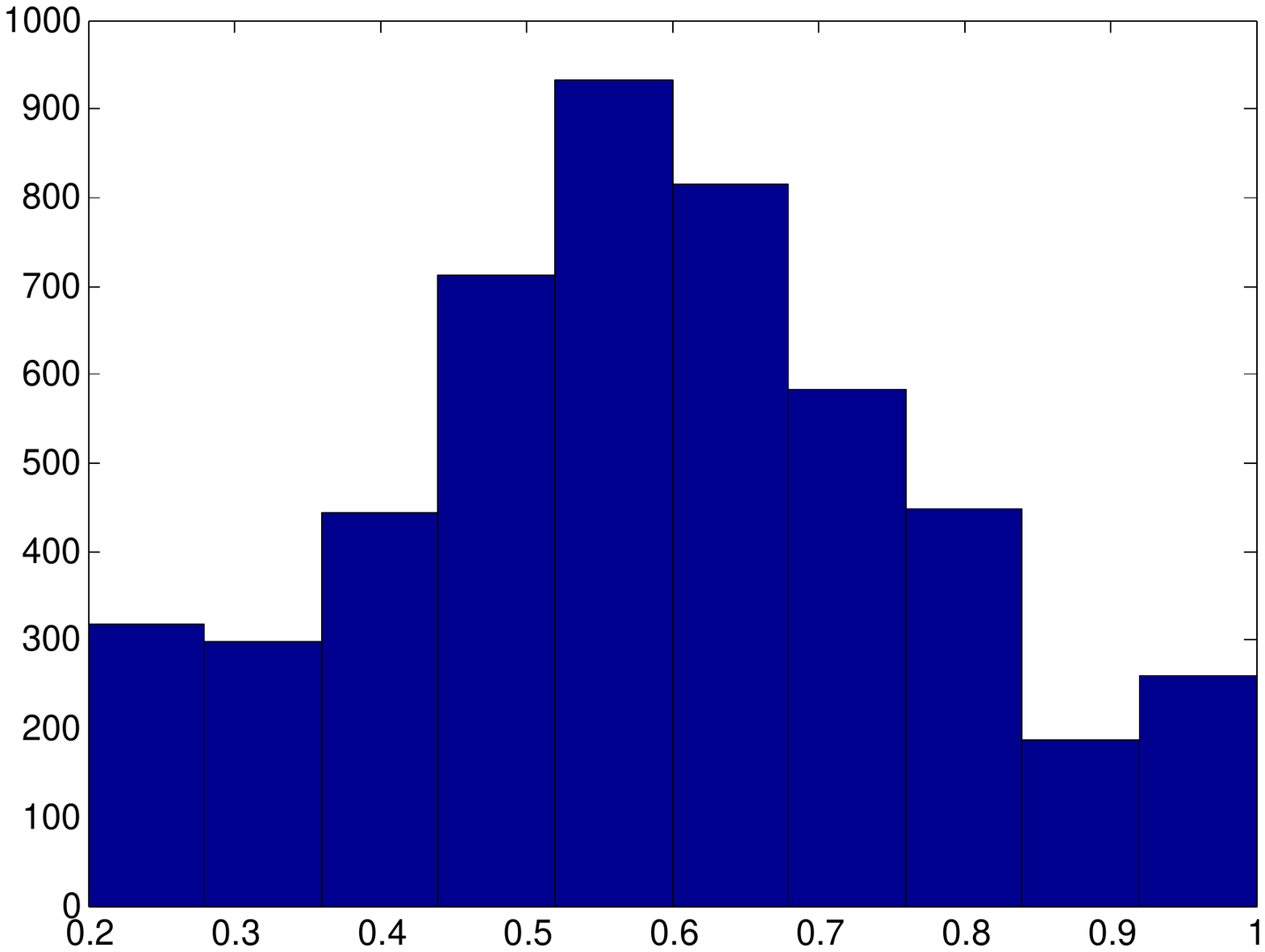}}

  \caption{ We plot a $10$ day sample path of an IVAPI policy using quadratic basis functions on the $1^{st}$ benchmark problem. Plot (a): The electricity price and resource level $R_{t}$. Plot (b): Histogram of the resource level.}
  \label{fig:Problem1}
\end{center}
\end{figure}

Figure \ref{fig:Problem1} depicts a sample path of an IVAPI policy for the $1^{st}$ benchmark problem in Table \ref{tab:benchmark}. The storage device is charged when electricity prices are low and discharged when electricity prices are high. We also note that the storage device fully discharges (below $20\%$) relatively infrequently.

\begin{figure}[htb]
  \begin{center}
  \includegraphics[width=6in,viewport=0 200 650 650]{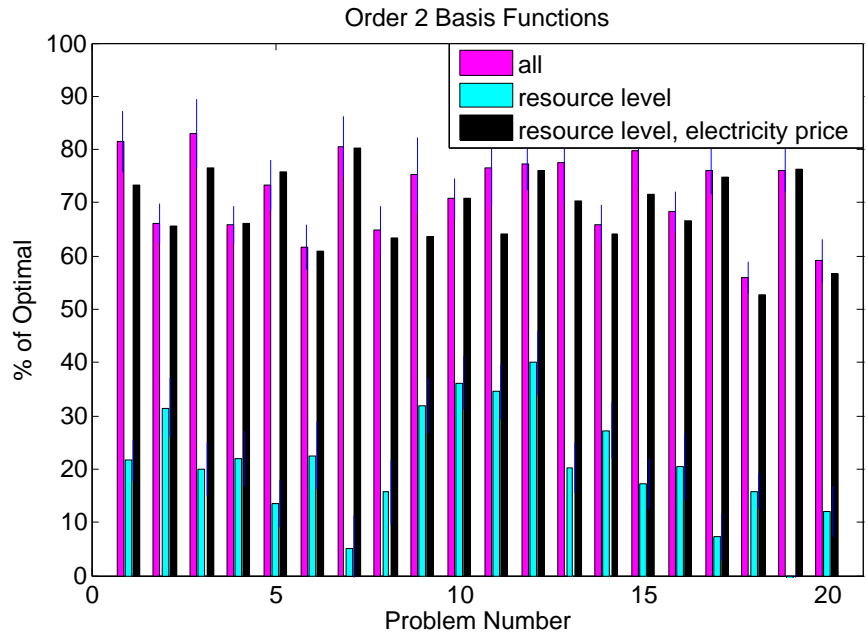}
  \end{center}
  \caption{The percentage of optimal along with 95\% confidence intervals for the average percent of optimal for the IVAPI policy, derived using quadratic basis functions, when only certain dimensions of the post-state are included in the post-state value function approximation.}
  \label{fig:leaveStatesOut}
\end{figure}

In order to reduce the number of basis functions in the algorithms, one may consider smaller dimensions for the post-decision state when constructing the value function approximation. Figure \ref{fig:leaveStatesOut} shows the results using three value function approximations: (1) all three state variables ($R_{t}$, $E_{t}$, and $P_{t}$), (2) only $R_{t}$, and (3) the resource and electricity prices $R_{t}$ and $P_{t}$. It is observed that using the resource level alone for the domain of the post-decision value function results in quite poor performances for most benchmark problems. Using both $R_{t}$ and $P_{t}$ appears to do fairly well overall, although using all of the state variable dimensions yields the best results. For certain problems, it may actually be advantageous to leave variables out of the value function approximation to simplify the policy search process.

\subsection{Benchmark Problems with Continuous State Spaces}\label{sec:RealisticProblem}
In this section we consider a set of $10$ problems with continuous state spaces, continuous actions, and the state transitions corresponding to the problem discussed in Section \ref{sec:Models}. These problems are listed in Table \ref{tab:benchmarkCont}. The electricity prices and loads are now time-dependent and stochastic for Problems 1-3.  Problems 4-10 are continuous steady-state problems. We compare both API algorithms and the myopic policy, although an optimal solution will no longer be available.

\begin{table}[htb]
\begin{center}
\begin{tabular}{|c|c|c|c|c|c|c|c|c|c|c|}
\hline
               \multicolumn{2}{|c|}{Problem}             & \multicolumn{5}{c|}{Number of Discretization Levels}& \multicolumn{4}{c|}{Parameters} \\
\hline

Number & Type & Time & Resource & Price & Load  & Wind   & Wind & Storage & RTE & Charge Rate \\
\hline
1    & Full & 96   & Cont.    & Cont. & Cont. & Cont.  & 0.1  & 2.5     & .81 & C/10 \\
\hline
2    & Full & 96   & Cont.    & Cont. & Cont. & Cont.  & 0.1  & 5.0     & .81 & C/10 \\
\hline
3    & BA   & 96   & Cont.    & Cont. & 1     & 1      & -    & -       & .81 & C/10 \\
\hline
4    & Full & 1    & Cont.    & Cont. & Cont. & Cont.  & 0.1  & 5.0     & .81 & C/10 \\
\hline
5    & Full & 1    & Cont.    & Cont. & Cont. & Cont.  & 0.1  & 2.5     & .81 & C/1 \\
\hline
6    & Full & 1    & Cont.    & Cont. & Cont. & Cont.  & 0.1  & 2.5     & .70 & C/1 \\
\hline
7    & BA   & 1    & Cont.    & Cont. & 1     & 1      & -    & -       & .81 & C/10 \\
\hline
8    & Full & 1    & Cont.    & Cont. & Cont. & Cont.  & 0.1  & 5.0     & .81 & C/1 \\
\hline
9    & Full & 1    & Cont.    & Cont. & Cont. & Cont.  & 0.1  & 5.0     & .70 & C/1 \\
\hline
10   & Full & 1    & Cont.    & Cont. & Cont. & Cont.  & 0.2  & 2.5     & .81 & C/1 \\
\hline
\end{tabular}
\caption{Parameter settings for problems with continuous states. Problems 1-3 have time-dependent stochastic loads and electricity prices.  Problems 4-10 are steady-state. \label{tab:benchmarkCont}}
\end{center}
\end{table}

\begin{figure}[htb]
\begin{center}
    \subfigure[]{\label{fig:AllContinuousProblems}\includegraphics[width=3in,viewport=50 200 550 550]{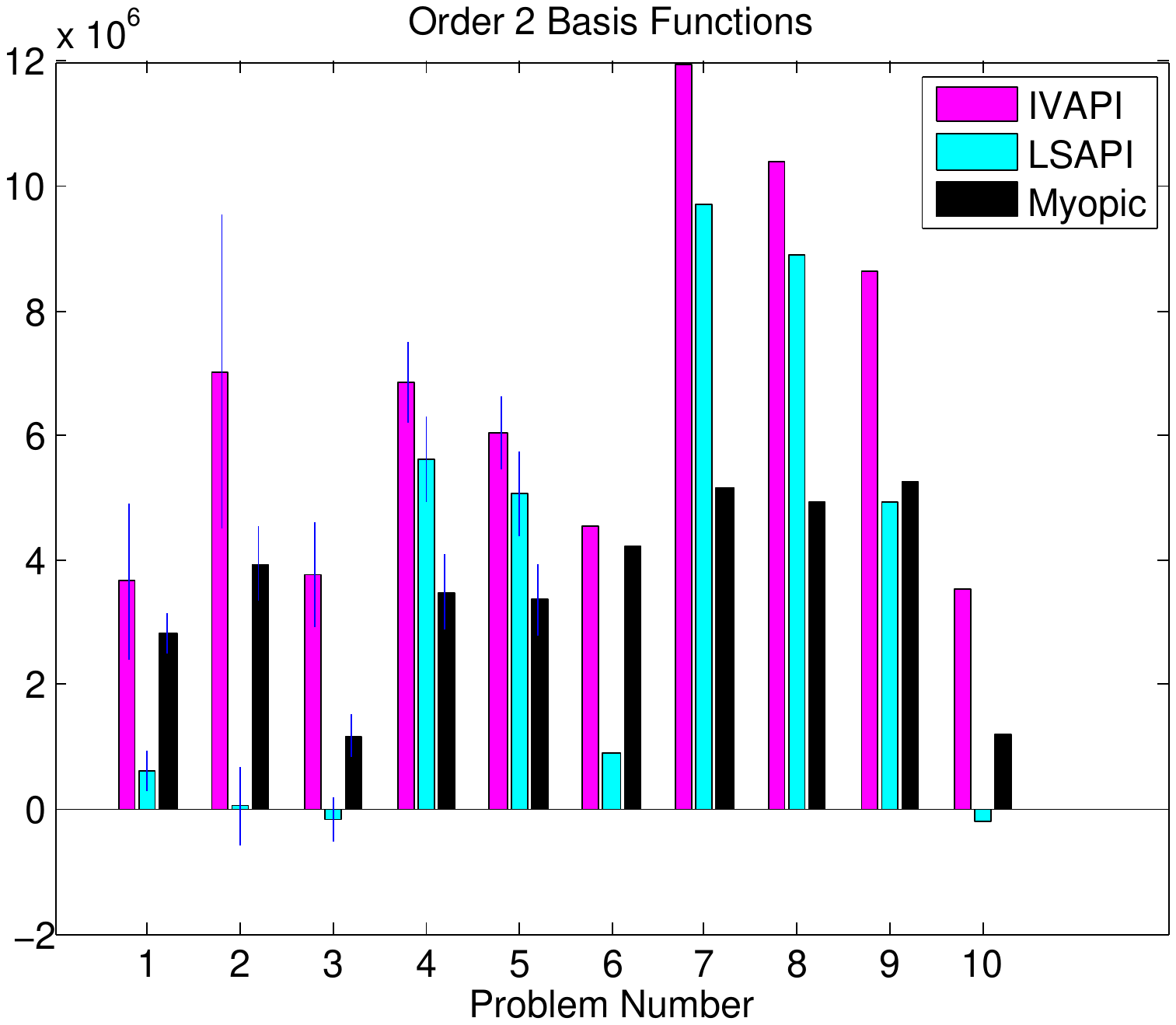}}
    \subfigure[] {\label{fig:continuousProblem5formatted}\includegraphics[width=3in,viewport=50 200 550 550]{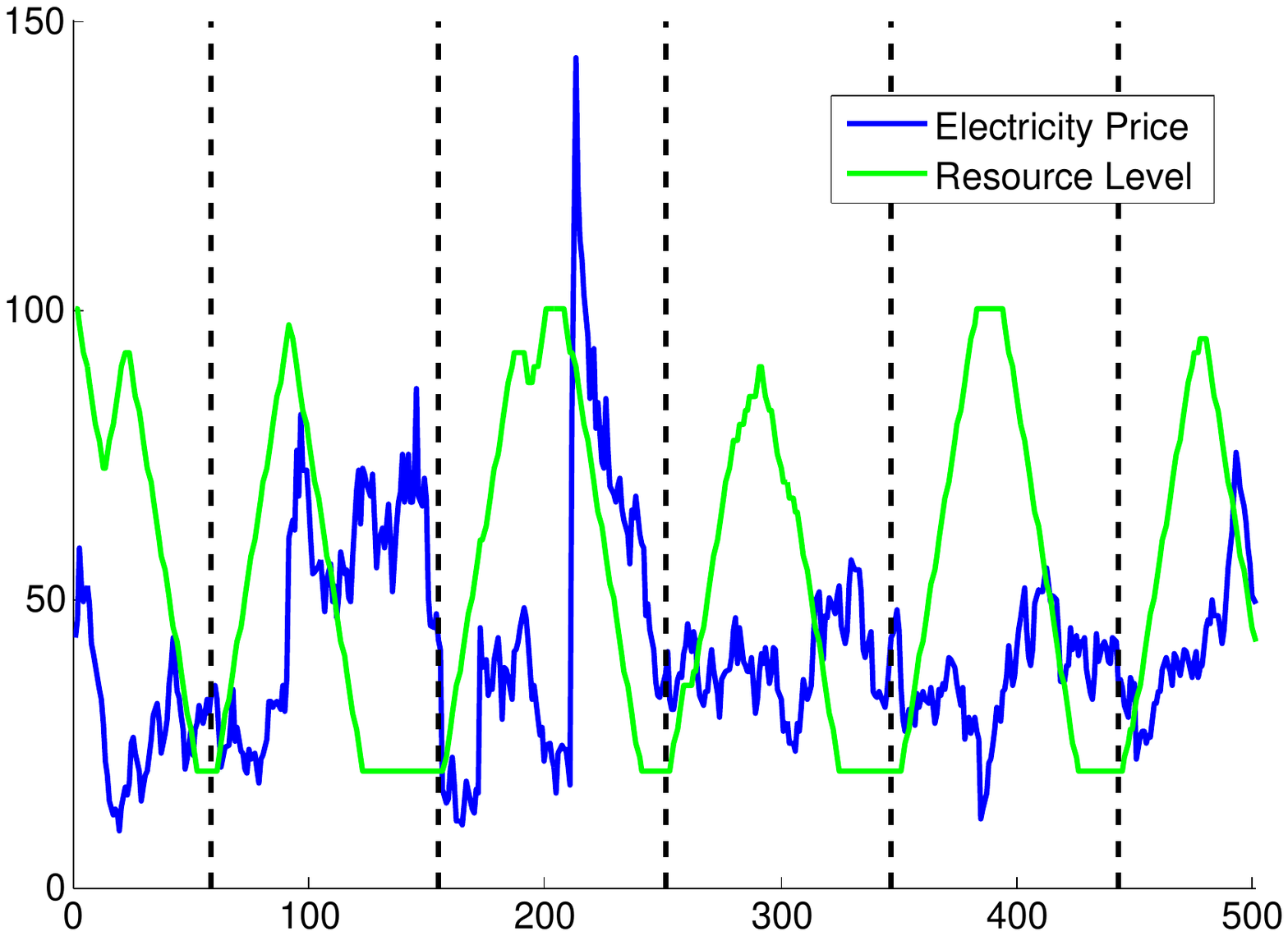}}
  \caption{Plot (a): Average objective of the myopic, IVAPI, and LSAPI policies on the continuous problems described in Table \ref{tab:benchmarkCont}. Plot (b): A sample path of the electricity spot price and resource level. The vertical lines correspond to midnight.}
\end{center}
\end{figure}

Figure \ref{fig:AllContinuousProblems} shows that for all $10$ test problems described in Table \ref{tab:benchmarkCont}, IVAPI prominently outperform LSAPI and the myopic policy. We see again that the use of instrumental variables in the API method adds value. Although all the dimensions of the state variable and action space are difficult to visualize, in Figure \ref{fig:continuousProblem5formatted} we use an IVAPI policy to show the electricity price and $R_{t}$ on one particular sample path. This figure shows that the policy tends to start charging the battery at night when electricity prices are low and then discharges the storage device throughout the day when electricity prices are higher.

The results in this section suggests that IVAPI is a promising and scalable algorithm, suitable for problems where the states are continuous and the transition probabilities are unknown.

\section{Conclusions}\label{Sec:conclusions}
In this paper, we have studied four variants of approximate policy iteration (API) methods, based on Bellman error minimization policy evaluation. We have considered least-squares Bellman error minimization, Bellman error minimization using instrumental variables, least-squares projected Bellman error minimization, and projected Bellman error minimization using instrumental variables. We have established the equivalence between Bellman error minimization using instrumental variables and projected Bellman error minimization methods. This result unifies algorithmic strategies that were previously thought to be fundamentally different.

The approximate policy iteration methods were then evaluated numerically. Motivated by control problems in energy storage applications, we created a series of benchmark problems to compare the different algorithmic strategies. The performance of a policy search method was also evaluated on those
benchmark problems.

We can draw several interesting conclusions from our numerical work. Bellman error minimization using instrumental variables appears to improve significantly over the basic LSAPI method, but otherwise did not work well when compared to the optimal benchmark, producing results that ranged between 60 and 80 percent of optimal.  Direct policy search performed much better, with results averaging over 90 percent of optimal.  Given that this is a problem that is ideally suited to LSAPI, it calls into question whether this is a method that can be counted on to produce good results.  The difficulty can probably be traced to the need to use pure exploration when sampling the state space, but given that the value function is relatively smooth and well approximated using a quadratic, it seems a bit surprising that the use of pure exploration would produce errors that are this large.

This research suggests that there are clear advantages to using direct policy search, possibly in conjunction with approximate policy iteration. The challenge is that in its derivative-free form, policy search does not scale easily with the dimension of the parameter space. This may be a major limitation in time-dependent applications, where we may need to estimate a different set of parameters for each time period.
\begin{APPENDICES}
\section{The Instrumental Variable Method}\label{app:IV}
The instrumental variable method (IVM) is a well known technique for dealing with errors in the explanatory variables of a regression problem, and provides a way to obtain consistent parameter estimates, see e.g., \citep{Du54,BrBa96,IV:Methods,KeSt61,IV:book:2005}. Consider the linear model in the matrix form:
\begin{eqnarray}
\label{eq:IVstructuralRelation}&&Y = X \beta, 
\end{eqnarray}
where $Y$ is a $n \times 1$ vector of response variables, $X$ is a $n \times k$ matrix of explanatory variables, and $\beta$ is a $k \times 1$ vector of weights. Let $X'$ and $Y'$ be observable values of the true values $X$ and $Y$. Following similar notation as in \citep{Du54}, denote the errors in the observed values of $X$ and $Y$ by $X''$  and $Y''$, respectively.
\begin{eqnarray}
\label{eq:IVerrorInX}&&X' = X + X'', \\
\label{eq:IVerrorInY}&&Y' = Y + Y''.
\end{eqnarray}
Hence, the linear regression model can now be written as $Y' = X \beta + Y''$. The least-squares method is probably the most widely used technique to estimate $\beta$. \cite{KeSt61} and \cite{Du54} show that least squares regression encounters problems with the model given by equations \eqref{eq:IVstructuralRelation}, \eqref{eq:IVerrorInX}, \eqref{eq:IVerrorInY}, since the observation of $X$ is typically correlated with the error in $X$. If $\beta$ is a scalar ($k=1$), this is easy to show.  Starting with the least squares estimate of $\beta_1$,
\begin{eqnarray}
\hat{\beta_1}=((X^{'})^{\top} X^{'})^{-1} (X^{'})^{\top}  Y^{'} &=& \frac{\sum_{i=1}^n X^{'}_{i1} Y^{'}_i}{\sum_{i=1}^n (X^{'}_{i1})^2} \nonumber\\
&=& \frac{\sum_{i=1}^n X^{'}_{i1}(X_{i1} \beta_1 + Y^{''}_i)}{\sum_{i=1}^n (X^{'}_{i1})^2} \label{eq:IVbiasLS1}\\
&=& \frac{\sum_{i=1}^n X^{'}_{i1} \left((X^{'}_{i1}- X^{''}_{i1}) \beta_1 + Y^{''}_i \right)}{\sum_{i=1}^n (X^{'}_{i1})^2} \label{eq:IVbiasLS2}\\
&=& \beta_1 - \beta_1 \frac{\sum_{i=1}^n X^{'}_{i1} X^{''}_{i1}}{\sum_{i=1}^n (X^{'}_{i1})^2} + \frac{\sum_{i=1}^n X^{'}_{i1} Y^{''}_i}{\sum_{i=1}^n (X^{'}_{i1})^2}, \label{eq:IVbiasLS3}
\end{eqnarray}
where equations \eqref{eq:IVstructuralRelation} and \eqref{eq:IVerrorInY} were used in equation \eqref{eq:IVbiasLS1}, and equation \eqref{eq:IVbiasLS2} follows \eqref{eq:IVerrorInX}. By taking the limit as $n\to \infty$ and assuming $\text{Cov}[X^{'}_{i1}, Y^{''}_{i}] =0$,  equation \eqref{eq:IVbiasLS3} becomes
\begin{eqnarray}
\lim_{n \to \infty} \hat{\beta_1} &=& \beta_1 - \beta_1 \lim_{n \to \infty}\left(\frac{\sum_{i=1}^n X^{'}_{i1} X^{''}_{i1}}{\sum_{i=1}^n (X^{'}_{i1})^2}\right) +  \lim_{n \to \infty} \left(\frac{\sum_{i=1}^n X^{'}_{i1} Y^{''}_i}{\sum_{i=1}^n (X^{'}_{i1})^2}\right) \nonumber\\
&=&\beta_1 - \beta_1 \lim_{n \to \infty}\left(\frac{\sum_{i=1}^n X^{'}_{i1} X^{''}_{i1}}{\sum_{i=1}^n (X^{'}_{i1})^2}\right) +  \frac{\text{Cov}[X^{'}_{i1}, Y^{''}_{i}]}{\mathbb{E}[(X^{'}_{i1})^2]} \nonumber\\
&=&\beta_1 - \beta_1 \lim_{n \to \infty}\left(\frac{\sum_{i=1}^n X^{'}_{i1} X^{''}_{i1}}{\sum_{i=1}^n (X^{'}_{i1})^2}\right). \label{eq:IVbiasLS5}
\end{eqnarray}
For many problems, $X^{'}_{i1}$ and $X^{''}_{i1}$ are positively correlated. Whence equation \eqref{eq:IVbiasLS5} implies that the least squares estimate $\hat{\beta_1}$ for $\beta_1$ is inconsistent. However, a consistent estimate can be obtained when an instrumental variable is available.

%


As it is discussed in the following section, a properly chosen instrumental variable can yield a consistent estimator for $\beta$. The following assumptions are made on the noise in $X$ and $Y$:
\begin{assumption}\label{ass:IVmeanYnoiseZero}$\mathbb{E}[Y_i^{''}]=0, \qquad  \text{for every  }i=1,\cdots,n$.\end{assumption}
\begin{assumption}\label{ass:IVmeanXnoiseZero}$\mathbb{E}[X_{ij}^{''}]=0, \qquad \text{for every  }i=1,\cdots,n,\text{ and }j=1,\cdots,k.$\end{assumption}

Unlike a standard linear regression problem in which a functional relationship to relate mean of the dependant variable to the regressor variable is considered, the IV method studies equation \eqref{eq:IVstructuralRelation} as a structural relation between the two random variables $X$ and $Y$, e.g., see \cite{Du54, KeSt61}.

Assume that an instrumental variable, $Z_j$, exists such that it is correlated with the true $X_j$ but uncorrelated with the errors in the observations of $X$ and $Y$.
\begin{assumption}
\label{ass:IVuncorrelatedWithYnoise} $\text{Cov}[Z_{ij}, Y^{''}_i]=0, \qquad \text{for every  }~i=1,\cdots,n,\text{ and }~j=1,\cdots,k.$
\end{assumption}
\begin{assumption}
\label{ass:IVuncorrelatedWithXnoise} $\text{Cov}[Z_{ij}, X^{''}_{il}]=0, \qquad \text{for every  }~i=1,\cdots,n,\text{ and }~j,l=1,\cdots,k.$
\end{assumption}
\begin{assumption}
\label{ass:IVsigmaFullRank} The $k \times k$ matrix $\Sigma$ has full rank $k$, where $X_l$ indicate the $l$'th column of $X$ and $\Sigma_{jl} = \text{Cov}[Z_j, X_l]$.
\end{assumption}
\begin{assumption}
\label{ass:IVuncorrelatedWithYnoiseLimit}$\lim_{n \to \infty}\frac{1}{n} \sum_{i=1}^n Z_{ij} Y^{''}_i = 0, \qquad \text{for every  }~j=1,\cdots,k.$
\end{assumption}
\begin{assumption}
\label{ass:IVuncorrelatedWithXnoiseLimit}$\lim_{n \to \infty}\frac{1}{n} \sum_{i=1}^n Z_{ij} X^{''}_{il} = 0, \quad \text{for every }~j,l=1,\cdots,k.$
\end{assumption}
\begin{assumption}
\label{ass:IVsigmaFullRankLimit}$\lim_{n \to \infty}\frac{1}{n} \sum_{i=1}^n Z_{ij} X_{il} = \text{Cov}[Z_j, X_l], \quad \text{for every }~j,l=1,\cdots,k.$
\end{assumption}

Assumptions \ref{ass:IVuncorrelatedWithYnoiseLimit}, \ref{ass:IVuncorrelatedWithXnoiseLimit}, and \ref{ass:IVsigmaFullRankLimit} do not follow trivially from assumptions \ref{ass:IVuncorrelatedWithYnoise}, \ref{ass:IVuncorrelatedWithXnoise}, \ref{ass:IVsigmaFullRank} without additional assumptions such as independence across the $n$ observations (because the law of large numbers does not apply). The method of instrumental variables defines the estimator $\hat{\beta}$ as the solution to
\begin{equation}Z^{\top} X' \hat{\beta} = Z^{\top} Y', \label{eq:instrumentalvariable1}\end{equation}
where $Z$ is a $n \times k$ matrix.  Note that $\hat{\beta}$ is uniquely defined when $Z^{\top} X^{'}$ has full rank $k$.

Below we extend the consistency proof from \cite{KeSt61} to use multiple instrumental variables ($k>1$).

\begin{proposition}
For the model given by equations \eqref{eq:IVstructuralRelation}, \eqref{eq:IVerrorInX}, \eqref{eq:IVerrorInY}, if assumptions \ref{ass:IVmeanYnoiseZero}, \ref{ass:IVmeanXnoiseZero}, \ref{ass:IVuncorrelatedWithYnoise}, \ref{ass:IVuncorrelatedWithXnoise}, \ref{ass:IVsigmaFullRank}, \ref{ass:IVuncorrelatedWithYnoiseLimit}, \ref{ass:IVuncorrelatedWithXnoiseLimit}, \ref{ass:IVsigmaFullRankLimit} hold, then $\hat{\beta} = [Z^{\top} X']^{-1} Z^{\top} Y'$  is a consistent estimator of $\beta$.
\label{thm:IVconsistent}
~\\
\textbf{Proof:} Using equation \eqref{eq:instrumentalvariable1} and equations \eqref{eq:IVstructuralRelation}, \eqref{eq:IVerrorInX}, and \eqref{eq:IVerrorInY}, we get
\begin{eqnarray}
(Z^{\top} X +  Z^{\top} X^{''}) \hat{\beta} =Z^{\top} (X + X^{''}) \hat{\beta}=Z^{\top} (X \beta + Y'')= Z^{\top} X \beta + Z^{\top} Y''.\label{eq:IVconsistency2}
\end{eqnarray}
Taking the limit as $n$ goes to infinity yields,
\begin{eqnarray}
\nonumber &&\lim_{n \to \infty}\frac{1}{n} \left(Z^{\top} X +  Z^{\top} X^{''}\right) \hat{\beta}=\lim_{n \to \infty} \frac{1}{n} \left( Z^{\top} X \beta + Z^{\top} Y'' \right).
\end{eqnarray}
By applying assumptions \ref{ass:IVuncorrelatedWithYnoise}  and \ref{ass:IVuncorrelatedWithYnoiseLimit} which imply $\lim_{n \to \infty}\frac{1}{n} Z^{\top} Y^{''} = \vec{0}$ and assumptions \ref{ass:IVuncorrelatedWithXnoise} and \ref{ass:IVuncorrelatedWithXnoiseLimit} which imply $\lim_{n \to \infty}\frac{1}{n} Z^{\top} X^{''} = \mathbf{0}$, and Slutsky's theorem when taking the limit, we arrive at:
\begin{eqnarray}
\nonumber&&\lim_{n \to \infty}\frac{1}{n} Z^{\top} X \hat{\beta}=\lim_{n \to \infty} \frac{1}{n} Z^{\top} X \beta.
\end{eqnarray}
Since the sample covariances converge in probability to the true covariances by Assumption \ref{ass:IVsigmaFullRankLimit}, we get $\Sigma \left(\lim_{n \to \infty} \hat{\beta}\right)=\Sigma \beta$. Now, Assumption \ref{ass:IVsigmaFullRank} which ensures that $\hat{\beta}$ is unique yields $\lim_{n \to \infty} \hat{\beta}=\beta$. \qed
\end{proposition}
\end{APPENDICES}






%
%
%




\end{document}